\documentclass[a4paper, 11 pt, twoside]{amsart}
\usepackage[left = 3.5cm, right = 3.5cm, headsep = 6mm,
footskip = 10mm, top = 35mm, bottom = 35mm, footnotesep=5mm, headheight =
2cm]{geometry}

\usepackage{multicol}

\usepackage{makecell}
\usepackage{mathtools,leftindex,tensor}
\usepackage[version=4]{mhchem}
\usepackage{longtable}
\usepackage{booktabs}
\usepackage{tabularx}

\RequirePackage{amsmath}
\RequirePackage{amssymb, amsthm, amsfonts}
\RequirePackage[T1]{fontenc}
\usepackage{comment}

\RequirePackage[british]{babel}

\RequirePackage{booktabs}
\RequirePackage{multirow}

\RequirePackage{enumitem}

\RequirePackage[usenames, dvipsnames, pdftex]{xcolor}
\RequirePackage{tikz, tikz-cd, tcolorbox}

\usetikzlibrary{shapes,backgrounds}
\usetikzlibrary{arrows, positioning, intersections}
\usetikzlibrary {through, chains}

\tikzset{every picture/.style={>=Straight Barb}} 
\tikzset{commutative diagrams/arrow style=tikz}

\newtheoremstyle{Lehn-it}
{}
{}
{\itshape}
{}
{\bfseries}
{$\;$\textmd{---}}
{ }
{}

\newtheoremstyle{Lehn-up}
{}
{}
{\upshape}
{}
{\bfseries}
{$\;$\textmd{---}}
{ }
{}

\newtheoremstyle{up-list}
{}
{}
{\upshape}
{}
{\bfseries}
{ }
{ }
{}

\newtheoremstyle{Lehn-Bemerkung}
{}
{}
{}
{}
{\itshape}
{$\;$\textmd{---}}
{ }
{}

\newtheoremstyle{citing}
{}
{}
{\itshape}
{}
{\bfseries}
{$\;$\textmd{---}}
{.5em}
{\thmnote{#3}}

\usepackage{chngcntr}
\usepackage{apptools}

\AtAppendix{\numberwithin{equation}{chapter}}
\AtAppendix{\numberwithin{equation}{section}}
\AtAppendix{\numberwithin{equation}{subsection}}
\AtAppendix{\numberwithin{equation}{subsubsection}}

\numberwithin{equation}{section}

\newenvironment{customthm}[1]
{\renewcommand\theinnercustomthm{#1}\innercustomthm}
{\endinnercustomthm}

\newenvironment{customprop}[1]
{\renewcommand\theinnercustomprop{#1}\innercustomprop}
{\endinnercustomprop}

{\theoremstyle{Lehn-it}

	\newtheorem{lem}[equation]{Lemma}
	
	\newtheorem{prop}[equation]{Proposition}

	\newtheorem{innercustomthm}{Theorem}
	\newtheorem{innercustomprop}{Proposition}
	
}
{\theoremstyle{Lehn-up}
	
	\newtheorem{exam}[equation]{Example}
	
}
{\theoremstyle{Lehn-Bemerkung}
	\newtheorem{rem}[equation]{Remark}
	
}

{\theoremstyle{up-list}
	
}

{\theoremstyle{citing}
	}

\RequirePackage{tikz-cd}
\RequirePackage{amsmath, amssymb, amsthm, amsfonts}

\newcommand{\del}{\partial}
\newcommand{\dbar}{\overline\del}

\DeclareMathOperator{\id}{\mathsf{id}}

\DeclareMathOperator{\End}{\mathsf{End}}
\DeclareMathOperator{\Der}{\mathsf{Der}}

\newcommand{\Wedge}{\textstyle\bigwedge}

\DeclareMathOperator{\Aut}{\mathsf{Aut}}
\DeclareMathOperator{\Out}{\mathsf{Out}}

\DeclareMathOperator{\im}{\mathsf{im}}

\DeclareMathOperator{\GL}{\mathsf{GL}}
\DeclareMathOperator{\SL}{\mathsf{SL}}

\DeclareMathOperator{\Kur}{Kur}

\DeclareMathOperator{\Diff}{Diff}

\DeclareMathOperator{\Gr}{Gr}
\DeclareMathOperator{\Sp}{\mathsf{Sp}}
\DeclareMathOperator{\Inn}{\mathsf{Inn}}

\DeclareMathOperator{\Ad}{Ad}

\renewcommand{\epsilon}{\varepsilon}

\DeclareFontFamily{OT1}{rsfs}{}
\DeclareFontShape{OT1}{rsfs}{n}{it}{<-> rsfs10}{}
\DeclareMathAlphabet{\curly}{OT1}{rsfs}{n}{it}

\newcommand{\kc}{{\mathcal C}}

\newcommand{\km}{{\mathcal M}}

\newcommand{\kp}{{\mathcal P}}

\newcommand{\kt}{{\mathcal T}}

\newcommand{\kv}{{\mathcal V}}

\newcommand{\kz}{{\mathcal Z}}

\newcommand{\IC}{{\mathbb C}}

\newcommand{\IN}{{\mathbb N}}

\newcommand{\IP}{{\mathbb P}}
\newcommand{\IQ}{{\mathbb Q}}
\newcommand{\IR}{{\mathbb R}}

\newcommand{\IZ}{{\mathbb Z}}
\newcommand{\gotha}{{\mathfrak a}}

\newcommand{\frg}{{\mathfrak g}}
\newcommand{\frh}{{\mathfrak h}}

\newcommand{\gothu}{{\mathfrak u}}

\DeclareMathOperator{\ad}{\mathsf{ad}}

\usepackage[unicode,bookmarks, pdftex]{hyperref}
\hypersetup{colorlinks=true, citecolor=NavyBlue ,linkcolor=NavyBlue,urlcolor=Orange, pdfpagemode=UseNone, breaklinks=true}

\usepackage{enumitem}
\newlist{prooflist}{description}{1}
\setlist[prooflist]{font=\normalfont \itshape, labelindent = \parindent, leftmargin = 0pt}

\usepackage{comment}

\makeatletter
\g@addto@macro\@floatboxreset\centering
\makeatother

\usepackage[backend=biber, style=alphabetic, url= false, doi= false, maxbibnames= 99, isbn= false]{biblatex}
\usepackage{csquotes}

\makeatletter
\def\ps@fixedhead{%
	\ps@headings
	\def\leftmark{\MakeUppercase{Konstantin Wehler}}
	\def\rightmark{\MakeUppercase{Classification of 3-dimensional complex nilmanifolds}}
}
\makeatother

\title[Classification of 3-dimensional complex nilmanifolds]{Classification of 3-dimensional complex nilmanifolds -- Teichmüller and moduli spaces}

\author{Konstantin Wehler}
\address{Konstantin Wehler\\FB 12/Mathematik und Informatik\\
	Philipps-Universit\"at Marburg\\
	Hans-Meerwein-Str. 6\\
	35032 Marburg\\
	Germany}
\email{konstantin.wehler@uni-marburg.de}

\begin{document}

	\begin{abstract}
		We classify complex nilmanifolds of dimension three by describing the Teichmüller and moduli spaces of left-invariant complex structures for every real six-dimensional nilmanifold.
	\end{abstract}
	
	\maketitle

\section*{Introduction}

A nilmanifold is a compact quotient of a real, simply connected, nilpotent Lie group $G$ by a lattice $\Gamma\subset G$. Nilmanifolds have long been established as a reliable source of examples in complex and differential geometry. This can be attributed to the general phenomenon that the geometry of a nilmanifold $\Gamma \backslash G$ is encoded, to a large extent, by the Lie algebra of $G$. In addition, real nilpotent Lie algebras have been completely classified up to dimension seven, and there are finitely many isomorphism classes up to dimension six \cite{Magnin86}. Therefore, the class of real six-dimensional nilmanifolds, in particular, serves as a common testing ground for the study of various geometric structures \cite{Latorre14, Angella18, Popovici18}, which often relies on the classification of complex structures on six-dimensional nilpotent Lie algebras~\cite{ceballos2016invariant}.

 For a real nilpotent Lie algebra $\frg$ of even dimension, we denote by 
\[
\kc(\frg) = \{J \colon \frg \to \frg \, | \, J^2= -\id, \, \textup{$J$ integrable}\}
\]
the space of complex structures on $\frg$. In general, it is already a difficult problem to determine whether or not $\kc(\frg)$ is empty, but there is a complete picture up to dimension six.
 Apart from the abelian Lie algebras there is only one nilpotent Lie algebra with a complex structure up to dimension four, whereas Salamon proved in~\cite{salamon01} that there are 18 isomorphism classes of nilpotent Lie algebras admitting a complex structure in dimension six. Building on numerous earlier works \cite{Ugarte07, Andrada11, Ugarte14Non-nilpotent}, a detailed analysis of these complex structures was then carried out in \cite{ceballos2016invariant}. Concretely, they classified complex structures on any given six-dimensional nilpotent Lie algebra $\frg$ up to Lie algebra automorphism, and thereby described the quotient $\kc(\frg)/\Aut(\frg)$.

Let $M = \Gamma \backslash G$ be an even-dimensional nilmanifold with underlying Lie algebra $\frg$ and suppose $\kc(\frg)$ is non-empty. Then every complex structure $J \in \kc(\frg)$ descends to a complex structure on $M$, which we also denote by $J$, and the resulting compact complex manifold $(M,J)$ is called a complex nilmanifold. Since $M$ is uniquely determined by the pair $(\frg, \Gamma)$, we say that $(M,J)$ is a complex nilmanifold of type $(\frg, \Gamma)$. 
The goal of the present paper is to provide a classification of three-dimensional complex nilmanifolds by describing both the moduli space
\[
\km(\frg, \Gamma) = \left\{ \textup{complex nilmanifolds $(M,J)$}\right\}/ \cong
\] 
of biholomorphism classes of complex nilmanifolds of type $(\frg, \Gamma)$, and the so-called Teichmüller space
\[
\kt(\frg) =  \left\{ \textup{complex nilmanifolds $(M,J)$}\right\}/ \sim
\]
for every real nilmanifold of dimension six. Here, the equivalence relation $\sim$ identifies two complex nilmanifolds if they are biholomorphic via a diffeomorphism in the identity component of the diffeomorphism group of $M$. The omission of the lattice $\Gamma$ from the notation for the Teichmüller space is justified by the description of these spaces given in \cite{Wehler26Teich} (see Section \ref{sect: Teichmüller and moduli}). 

Since an automorphism of $\frg$ does not usually descend to a diffeomorphism of $M$, the spaces $\km(\frg,\Gamma)$ and $\kt(\frg)$ differ quite drastically from the moduli space $\kc(\frg)/\Aut(\frg)$ of complex structures on $\frg$. In addition to the case where $\frg$ is abelian, this is already illustrated by the only other four-dimensional nilpotent Lie algebra with a complex structure. In this case, $\km(\frg, \Gamma)$ and $\kt(\frg)$ can be described as complex two-dimensional manifolds \cite{Borcea84}, while $\kc(\frg)/\Aut(\frg)$ is just a point.

Following Ugarte's notation \cite{Ugarte07}, we denote the 18 six-dimensional nilpotent Lie algebras admitting a complex structure by $\frh_1, \dots, \frh_{16}, \frh_{19}^-$ and $\frh_{26}^+$, where $\frh_1 \cong \IR^6$ is just the abelian Lie algebra, and $\frh_5$ is the real Lie algebra underlying the Iwasawa manifold. If $\frg$ is one of these six-dimensional Lie algebras, then an associated simply connected Lie group $G$ always admits a lattice $\Gamma \subset G$ \cite{malcev51}. Moreover, the space $\kc(\frg)$ can be viewed as complex analytic subspace of the Grassmannian $\Gr(3, \frg_\IC)$, and thus the quotients $\km(\frg, \Gamma)$ and $\kt(\frg)$ are naturally topological spaces. 
In \cite{salamon01} Salamon provides an upper bound for the complex dimension of $\kc(\frg)$ for every six-dimensional nilpotent Lie algebra, which coincides with the actual dimension of $\kc(\frg)$ in all cases except $\frg = \frh_{15}$ and $\frg = \frh_{16}$ (Remark \ref{rem: dim C(g)}). 

We will determine those pairs $(\frg, \Gamma)$ for which $\km(\frg, \Gamma)$ or $\kt(\frg)$ can also be endowed with the structure of a complex space.
Apart from the abelian Lie algebra $\frh_1$, a complex structure was already established on $\kt(\frg)$ for $\frg \in \{\frh_2, \frh_3, \frh_4, \frh_6, \frh_8, \frh_9, \frh_{10}\}$ in \cite{Wehler26Teich} by classifying the associated complex nilmanifolds by the cohomology class of a trivialising section of their canonical bundle. More precisely, since the canonical bundle of every complex nilmanifold can be trivialised by a left-invariant section \cite{salamon01,BDV}, one can define a period map
\[
\kp \colon \kt(\frg) \to \IP H^3(\frg, \IC),
\]
which was shown to be injective in all of the above cases. The following theorem completes the description of $\kt(\frg)$ for all nilpotent Lie algebras of dimension six.

\begin{customthm}{A}\label{thm: teichmüller and period maps dim 6}
	Let $\frg$ be a real nilpotent Lie algebra of dimension six admitting complex structures. Then the period map $\kp \colon \kt(\frg) \to \IP H^3(\frg, \IC)$ is injective if and only if $\frg \not\in \{\frh_5,\frh_{13}, \frh_{14},\frh_{15},\frh_{16}\}$. Moreover, we have the following cases.
	\begin{enumerate}
		\item If $\frg \not\in \{ \frh_5, \frh_{15}, \frh_{16}\}$, then both $\kc(\frg)$ and the Teichmüller space $\kt(\frg)$ are complex manifolds.
		\begin{center}
			\begin{tabular}{l c}
				\toprule
				\textup{Lie algebra} & $\dim_\IC \kt(\frg)$\\
				\midrule
				$\frh_1$ & $9$\\
				$\frh_2,\frh_3,\frh_4$ & $4$\\
				$\frh_6,\frh_7$ & $5$\\
				$\frh_8$ & $6$\\
				$\frh_9,\frh_{10},\frh_{11},\frh_{12}, \frh_{13}, \frh_{14}$ & $3$\\
				$\frh_{19}^-,\frh_{26}^+$& $2$\\
				\bottomrule
			\end{tabular}
		\end{center}
		\item If $\frg \in \{\frh_5, \frh_{15}, \frh_{16}\}$, then the complex space $\kc(\frg)$ is not reduced, and the Teichmüller space $\kt(\frg)$ is non-Hausdorff for $\frg= \frh_5$ and $\frg = \frh_{15}$.
	\end{enumerate}
\end{customthm}

It turns out that the Hausdorff property of the Teichmüller space $\kt(\frg)$ has a geometric characterisation for any nilpotent Lie algebra $\frg$.

\begin{customprop}{B}\label{thm: h^0 constant}
	Let $\frg$ be a real nilpotent Lie algebra admitting complex structures, and let $G$ be an associated simply connected Lie group. If the Teichmüller space $\kt(\frg)$ is Hausdorff, then for any lattice $\Gamma\subset G$, the dimension of the automorphism group $\Aut(\Gamma\backslash G, J)$ is constant on the connected components of $\kc(\frg)$.
\end{customprop}

While Theorem \ref{thm: teichmüller and period maps dim 6} shows that the Teichmüller space has a nice geometric structure for most Lie algebras of dimension six, this is far from true for the moduli space. Consider the natural map
\[
\pi \colon \kt(\frg) \to \km(\frg, \Gamma),
\]
whose fibres are exactly the points in $\kt(\frg)$ corresponding to biholomorphic complex nilmanifolds. It was already proved by Kodaira and Spencer in \cite{Kodaira58} that the moduli space for complex tori of dimension $n$ is non-Hausdorff for $n \geq 2$. In fact, they showed that the fibre of $\pi$ over a generic complex torus is dense.
Up to dimension six, the Lie algebra $\frh_7$ is the only non-abelian Lie algebra for which $\pi$ possibly has dense fibres (Proposition \ref{prop: dense orbits}). Nevertheless, the moduli space usually contains non-closed points.

\begin{customthm}{C}\label{thm: coarse moduli space dim 6}
	Let $\frg$ be a real nilpotent Lie algebra of dimension six admitting complex structures, and let $G$ be an associated simply connected Lie group. Then the following holds for any lattice $\Gamma \subset G$.
	\begin{enumerate}
		\item If $\frg \not\in \{ \frh_{19}^-, \frh_{26}^+\}$, then the moduli space $\km(\frg, \Gamma)$ contains non-closed points, and is thus not Hausdorff. Moreover, if $\frg \in \{ \frh_1, \frh_6, \frh_7, \frh_8, \frh_9, \frh_{10}\}$, then the set of non-closed points in $\km(\frg, \Gamma)$ is dense.
		\item If $\frg \in \{\frh_{19}^-, \frh_{26}^+\}$, then the moduli space $\km(\frg, \Gamma)$ is a normal complex space of dimension $2$.
	\end{enumerate}
\end{customthm}

Complex nilmanifolds of type $(\frh_3, \Gamma)$ belong to the class of so-called Kodaira manifolds whose moduli space was already known to be non-Hausdorff \cite{grantcharov2004deformations}, with the exception of the moduli space for Kodaira surfaces which was shown to be a complex manifold by Borcea \cite{Borcea84}.
In the proof of Theorem \ref{thm: coarse moduli space dim 6}, we study the topology of the moduli spaces $\km(\frg, \Gamma)$ using various criteria invoking the arithmetic properties of the automorphism group $\Aut(\Gamma)$ (see Section \ref{subsect: moduli criteria}).

Among the six-dimensional nilpotent Lie algebras with a complex structure, $\frh_{19}^-$ and $\frh_{26}^+$ are the only ones with one-dimensional centre. Although there is no general classification of eight-dimensional nilpotent Lie algebras admitting a complex structure, there is a classification of those with one-dimensional centre \cite{Latorre23}. In light of Theorem \ref{thm: coarse moduli space dim 6}, one might hope to be able to construct nice moduli spaces for some Lie algebras in this class as well.

\subsection*{Acknowledgements}
This paper is a revised version of a chapter from my PhD thesis. I would like to express my gratitude to my advisers Sönke Rollenske and Nicolina Istrati for their constant support, and many valuable comments and discussions. This work was supported by the DFG through the grant RO 3734/4-1. 
After the completion of this paper, Claude was used to provide a Python script to verify the computations in Section \ref{sect: moduli dim six}.

\section{Preliminaries on complex nilmanifolds}

In this section, we recall some notation and well-known results about complex nilmanifolds, particularly in complex dimension three.

\subsection{Complex structures on nilpotent Lie algebras}\label{sect: complex structures on nilpotent Lie algebras}

Let $\frg$ be a real nilpotent Lie algebra of dimension $2n$. An almost complex structure $J \colon \frg \to \frg$ is said to be integrable if 
\begin{equation*}\label{eq: nijenhuis}
	[x,y]-[Jx,Jy]+J[Jx,y]+J[x,Jy] =0
\end{equation*}
for all $x,y \in \frg$. In this case, we say that $J$ is a complex structure on $\frg$. By taking the dual of the Lie bracket on $\frg$, we obtain a map $d \colon \frg^* \to \Wedge^2 \frg^*$, which extends to a map on the exterior algebra $\Wedge^\bullet \frg^*$. The Lie algebra $\frg$ is then determined by the real structure equations
\[
de^k = \sum_{i<j} c_{ij}^k e^i \wedge e^j
\]
for a given basis  $e^1, \dots, e^{2n}$ of $\frg^*$.
Following Salamon's notation we conveniently express these equations by the tuple $(de^1, \dots, de^{2n})$, and abbreviate $e^i \wedge e^j$ to $ij$. For instance, $(0,0,0,12)$ describes the four-dimensional nilpotent Lie algebra whose only non-zero differential is $de^4 = e^1 \wedge e^2$, or equivalently, if $e_1, \dots, e_4$ is the dual basis, then the only non-zero bracket is $[e_1,e_2]= -e_4$. While the Lie algebra $(0,0,0,12)$ is the only non-abelian nilpotent Lie algebra up to dimension four admitting a complex structure~\cite[Proposition~2.3]{salamon01}, Salamon proved in \cite{salamon01} that a complex structure exists on exactly the following 18 six-dimensional nilpotent Lie algebras:
\begin{enumerate}
	\begin{multicols}{2}
		\item[] $\mathfrak{h}_1 = (0,0,0,0,0,0)$
		\vspace{2pt}
		\item[] $\mathfrak{h}_2 = (0,0,0,0,12,34)$
		\vspace{2pt}
		\item[] $\mathfrak{h}_3 = (0,0,0,0,0,12+34)$
		\vspace{2pt}
		\item[] $\mathfrak{h}_4 = (0,0,0,0,12,14+23)$
		\vspace{2pt}
		\item[] $\mathfrak{h}_5 = (0,0,0,0,13+42,14+23)$
		\vspace{2pt}
		\item[] $\mathfrak{h}_6 = (0,0,0,0,12,13)$
		\vspace{2pt}
		\item[] $\mathfrak{h}_7 = (0,0,0,12,13,23)$
		\vspace{2pt}
		\item[] $\mathfrak{h}_8 = (0,0,0,0,0,12)$
		\vspace{2pt}
		\item[] $\mathfrak{h}_9 = (0,0,0,0,12,14+25)$\\
		\item[] $\mathfrak{h}_{10} = (0,0,0,12,13,14)$
		\vspace{2pt}
		\item[] $\mathfrak{h}_{11} = (0,0,0,12,13,14+23)$
		\vspace{2pt}
		\item[] $\mathfrak{h}_{12} = (0,0,0,12,13,24)$
		\vspace{2pt}
		\item[] $\mathfrak{h}_{13} = (0,0,0,12,13+14,24)$
		\vspace{2pt}
		\item[] $\mathfrak{h}_{14} = (0,0,0,12,14,13+42)$
		\vspace{2pt}
		\item[] $\mathfrak{h}_{15} = (0,0,0,12,13+42,14+23)$
		\vspace{2pt}
		\item[] $\mathfrak{h}_{16} = (0,0,0,12,14,24)$
		\vspace{2pt}
		\item[] $\mathfrak{h}_{19}^- = (0,0,0,12,23,14-35)$
		\vspace{2pt}
		\item[] $\mathfrak{h}_{26}^+ = (0,0,12,13,23,14+25)$
	\end{multicols}
\end{enumerate}
\vspace*{-15pt}

 After linearly extending an almost complex structure $J\colon \frg \to\frg$ to the complexification $\frg_\IC$, we obtain a decomposition $\frg_\IC = \frg^{1,0} \oplus \frg^{0,1}$ into the $\pm i$-eigenspaces. Similarly, we can decompose its dual $\mathfrak{g}_\mathbb{C}^* = \mathfrak{g}^{*1,0} \oplus \mathfrak{g}^{*0,1}$, and we denote $\Wedge^{p,q}\mathfrak{g}^* = \Wedge^p\mathfrak{g}^{*1,0} \otimes \Wedge^q\mathfrak{g}^{*0,1}$. Therefore the exterior powers of $\mathfrak{g}_\IC^*$ decompose as
\begin{equation*}
	\Wedge^k\mathfrak{g}_\mathbb{C}^* = \bigoplus_{p+q = k} \Wedge^{p,q}\mathfrak{g}^*,
\end{equation*}
and moreover $\overline{\Wedge^{p,q}\frg^*} = \Wedge^{q,p}\frg^*$. We denote the projection $\Wedge^\bullet \frg_\IC^* \to \Wedge^{p,q} \frg^*$ by $\pi^{p,q}$. The differential $d$ defined above extends to the complexified exterior algebra $\Wedge^\bullet\mathfrak{g}_\mathbb{C}^*$, and we define
\[
\dbar = \pi^{p,q+1} \circ d \colon \Wedge^{p,q}\frg^* \to \Wedge^{p,q+1}\frg^*
\]
and
\[
\partial = \pi^{p+1,q} \circ d \colon \Wedge^{p,q}\frg^* \to \Wedge^{p+1,q}\frg^*.
\]
The almost complex structure $J$ is then integrable if and only if $d = \partial +\dbar$. As for the real Lie algebra $\frg$, we will often describe a complex structure on $\frg$ in terms of forms. Suppose $J$ is integrable and $\omega^1, \dots, \omega^n$ is a basis of $\frg^{*1,0}$. Then we have 
\[
d\omega^k = \sum_{i<j} A_{ij}^k \omega^{i} \wedge \omega^j +\sum_{i,j} B_{ij}^k \omega^{i} \wedge \overline \omega^j,
\]
where $A_{ij}^k, B_{ij}^k\in \IC$. These equations are called the complex structure equations with respect to the basis $\omega^1, \dots, \omega^n$, and they determine the complex structure $J$ up to Lie algebra automorphism. Throughout the text, we will usually abbreviate $\omega^i \wedge \omega^j$ to $\omega^{ij}$ and $\omega^i \wedge \overline \omega^j$ to $\omega^{i\bar j}$. For example, every complex structure on $(0,0,0,12)$ admits a basis $\omega^1, \omega^2$ with complex structure equations $d\omega^1 =0$ and $d\omega^2 = \omega^{1\bar1}$.
 The possible complex structure equations for all six-dimensional Lie algebras in the list above were completely described in \cite[Tables~1 and 2]{ceballos2016invariant}.
 
 We will also make use of the following cohomology attached to a complex structure $J$ on $\frg$. Following \cite{rollenske09b}, the so-called Dolbeault cohomology of $J$ is the cohomology $H^{p,\bullet}(\frg, J)$ of the complex
 \[
 \begin{tikzcd}
 	0  \rar & \Wedge^{p,0}\frg^* \rar{\dbar} & \Wedge^{p,1}\frg^* \rar{\dbar} &   {\cdots}  \rar{\dbar} & \Wedge^{p,n-1}\frg^* \rar{\dbar} &\Wedge^{p,n}\frg^* \rar & 0.
 \end{tikzcd}
 \]
 The differential $\dbar$ on $\Wedge^{0,\bullet}\frg^*$ can also be extended to $\frg^{1,0} \otimes \Wedge^{0,\bullet}\frg^*$ by defining the map $\dbar \colon \frg^{1,0} \to \frg^{1,0} \otimes \frg^{*0,1}$ as $\dbar(x)(\overline y) = -[x,\overline y]^{1,0}$ for $x \in \frg^{1,0}$ and $\overline y \in \frg^{0,1}$. We denote the cohomology of the associated complex $(\frg^{1,0} \otimes \Wedge^{0,\bullet} \frg^*, \dbar)$ by $H^\bullet(\frg, \frg^{1,0})$.
\subsection{Complex nilmanifolds}

Let $\frg$ be a real nilpotent Lie algebra and let $G$ be an associated simply connected Lie group. A nilmanifold is a compact quotient $M= \Gamma \backslash G$, where $\Gamma\subset G$ is a lattice, that is, $\Gamma$ is a discrete co-compact subgroup of $G$. It was shown by Malcev that such a lattice exists if and only if $\frg^*$ admits a basis whose structure equations have rational coefficients \cite{malcev51}. Note that all of the six-dimensional nilpotent Lie algebras (with a complex structure) admit such a basis. 

If $M$ is even-dimensional and $J$ is a complex structure on $\frg$, then we can extend $J$ to a complex structure on $G$ by left-translations, and since $\Gamma$ acts on the left this complex structure descends to $M$. We also denote this complex structure by $J$, and we refer to complex structures on $M$ obtained in this way as left-invariant.

 A complex nilmanifold $X = (M,J)$ is a nilmanifold $M$ endowed with a left-invariant complex structure $J$. To specify the diffeomorphism type of $X$, we sometimes say that $X$ is a complex nilmanifold of type $(\frg, \Gamma)$ if its underlying smooth nilmanifold is the quotient of a simply connected Lie group with Lie algebra $\frg$ by a lattice $\Gamma$.

It is a well-known result due to Nomizu \cite{nomizu54} that the de Rham cohomology of a nilmanifold with underlying Lie algebra $\frg$ is computed by the complex $\left(\Wedge^\bullet \frg^*, d\right)$. More recently, it was shown in \cite{Has26} that the Dolbeault cohomology of any complex nilmanifold can also be computed entirely on its underlying Lie algebra. Concretely, if $X$ is an $n$-dimensional complex nilmanifold of type $(\frg, \Gamma)$ with complex structure $J$, then its Dolbeault cohomology $H^{p,q}_{\dbar}(X)$ is isomorphic to $H^{p,q}(\frg,J)$.
Similarly, the cohomology of $X$ with values in the sheaf of holomorphic vector fields $H^q(X, \Theta_X)$ is isomorphic to $H^q(\frg, \frg^{1,0})$ \cite[Corollary~2.5]{rollenske09b}.

\section{Teichmüller and moduli spaces}\label{sect: Teichmüller and moduli}

Following \cite{Wehler26Teich}, we recall some basic results about the Teichmüller and moduli spaces for complex nilmanifolds. A complex structure on a real nilpotent Lie algebra $\frg$ of dimension $2n$ is uniquely determined by a complex $n$-dimensional subspace $V \subset \frg_\IC$ satisfying $V \cap \overline V = 0$ and the integrability condition $[V,V] \subset V$. Therefore, we can identify the space of complex structures on $\frg$ with the subspace
\[
\kc(\frg) =\{ V \in \Gr(n, \frg_\IC) \, | \, V \cap \overline{V} =0,\, [V, V]\subset V\}
\]
of the Grassmannian $\Gr(n, \frg_\IC)$. In particular, $\kc(\frg)$ has the structure of a complex space, which can be arbitrarily singular \cite[Section~3]{Wehler26Teich}. 

For any nilmanifold $M = \Gamma \backslash G$ with underlying Lie algebra $\frg$, the space $\kc(\frg)$ contains all left-invariant complex structures on $M$. So we can describe the Teichmüller and moduli spaces as the quotients $\kt(\frg) = \kc(\frg)/\sim$ and $\km(\frg, \Gamma) = \kc(\frg)/\cong$. Here, the equivalence relation $\cong$ identifies two complex nilmanifolds $(M,J)$ and $(M,J')$ if there exists a biholomorphism $f \colon (M,J) \to (M,J')$, whereas $\sim$ requires $f$ to be in the identity component of the diffeomorphism group $\Diff(M)$ of $M$.

There is, however, a more accessible description of these two spaces, proved in~\cite{Wehler26Teich}. The space $\kc(\frg)$ admits several natural group actions. Since we take the quotient by $\Gamma$ from the left, every right-translation on $G$ descends to a diffeomorphism of $M$, which is in the identity component of $\Diff(M)$. Therefore, the group $G$ acts on the space $\kc(\frg)$, viewed as the space of left-invariant complex structures on $M$, by right-translations. Furthermore, any automorphism of the lattice $\Gamma$ extends uniquely to an automorphism of $G$ \cite[Corollary~2.8]{Vinberg00}, which then induces a diffeomorphism of $M$. In this way, the discrete group $\Aut(\Gamma)$ also acts on $\kc(\frg)$. This action further descends to the quotient $\kc(\frg)/G$, and it was shown in \cite[Corollary~4.3]{Wehler26Teich} that
\[
\kt(\frg) = \kc(\frg)/G
\]
and 
\[
\km(\frg, \Gamma) = \kt(\frg)/\Aut(\Gamma).
\]
We endow both of these spaces with the quotient topology of $\kc(\frg) \subset \Gr(n, \frg_\IC)$.

\begin{rem} 
	In the present text, we restrict our considerations to left-invariant complex structures on $M$, but one could also consider the (classical) Teichmüller space $\kt(M)= \kc(M)/\sim$, where $\sim$ is defined as above and $\kc(M)$ is the space of all complex structures on $M$.
	
	In general, the space $\kt(\frg) \subset \kt(M)$ is just an open subset, but $\kt(\frg)$ forms a union of connected components in $\kt(M)$ for several classes of Lie algebras, including all abelian Lie algebras \cite{Cat04DefLargeI}, and all six-dimensional Lie algebras, with the exception of $\frg = \frh_7$ \cite{rollenske09}.
\end{rem}

\begin{rem}\label{rem: orientation}
	As explained in \cite[Remark~2.1]{Meersseman19} or \cite[Remark~4.5]{Wehler26Teich}, after fixing an orientation on $M = \Gamma \backslash G$ the number of connected components of $\km(\frg, \Gamma)$ depends on the existence of an orientation-reversing diffeomorphism of $M$. We will see in Example \ref{exam: lattice dependence} that the existence of such a diffeomorphism, and thus the number of connected components of $\km(\frg, \Gamma)$, does depend on the choice of $\Gamma \subset G$.
\end{rem}

One of our main goals is to determine when the topological spaces $\kt(\frg)$ and $\km(\frg, \Gamma)$ can be endowed with a complex structure, which is natural in the sense that it is compatible with the charts given by Kuranishi theory.

 More precisely, if $X$ is a complex nilmanifold of type $(\frg, \Gamma)$, then there exists a semi-universal family of deformations of $X$ over a finite-dimensional complex space $\Kur(X)$, known as the Kuranishi space of $X$, whose Zariski tangent space has dimension $\dim H^1(X, \Theta_X)$ \cite{kuranishi62}. This family is called the Kuranishi family of $X$, and it contains all sufficiently small deformations of $X$. Moreover, every sufficiently small deformation of $X$ is again a complex nilmanifold~\cite{rollenske09b, Has26}. Therefore, we have a map 
\[
 \pi_X \colon \Kur(X) \to \kt(\frg),
\]
which is surjective onto the germ of $\kt(\frg)$ at the point corresponding to $X$, but it is generally not a homeomorphism onto its image \cite{catanese2011superficial}. For instance, if $X$ is the Iwasawa manifold, then $\pi_X$ has fibres of positive dimension. 

If, however, the map $\pi_X$ is injective and the Kuranishi family is universal for all complex nilmanifolds of type $(\frg, \Gamma)$, then $\pi_X$ is a local homeomorphism and we can glue the complex spaces $\Kur(X)$ to a global complex structure on $\kt(\frg)$ (see \cite[Remark~14]{catanese2011superficial} or \cite[Proposition~4.16]{Wehler26Teich}). If $\kt(\frg)$ is also Hausdorff, this complex structure makes $\kt(\frg)$ into a complex space, which descends to a complex structure on $\km(\frg, \Gamma)$ under some assumptions on the $\Aut(\Gamma)$-action (see Section \ref{subsect: moduli criteria}). In this case, $\kt(\frg)$ and $\km(\frg, \Gamma)$ can also be described as coarse moduli spaces for certain moduli functors \cite[Section~4]{Wehler26Teich}.

\section{Existence criteria}

In this section, we provide some existence and non-existence criteria for a complex structure on the Teichmüller and moduli spaces.

\subsection{The Teichmüller space}
Let $\frg$ be a real nilpotent Lie algebra of dimension $2n$ admitting complex structures. Our main tools to describe the Teichmüller space are the following period maps studied in \cite{Wehler26Teich}. The map
\[
\widetilde{\kp}\colon \kc(\frg) \to \IP H^n(\frg, \IC)
\]
sends $J \in \kc(\frg)$ to the point in $\IP H^n(\frg, \IC)$ corresponding to the one-dimensional subspace $H^{n,0}(\frg, J) \subset H^n(\frg, \IC)$, and thus it descends to a global period map
\[
\kp \colon \kt(\frg) \to \IP H^n(\frg, \IC).
\]
 In general, both of these maps are just continuous, but $\widetilde \kp$ is holomorphic after passing to the reduction of $\kc(\frg)$ \cite[Proposition~5.2]{Wehler26Teich}, while $\kt(\frg)$ does not necessarily admit any complex structure. As above, we can also define the so-called local period map 
 \[
 \kp_X \colon \Kur(X) \to \IP H^n(\frg, \IC)
 \]
for any complex nilmanifold $X$ with underlying Lie algebra $\frg$, yielding the commutative diagram 
\[
\begin{tikzcd}
	\Kur(X) \arrow[r, "\pi_X"] \arrow[rd,swap,"\kp_X"] & \kt(\frg) \arrow[d,"\kp" ]\\
	&  \IP H^n(\frg, \IC).
\end{tikzcd}
\]
Furthermore, if the local period map $\kp_X$ is injective, then $\pi_X$ is a local homeomorphism \cite[Proposition~15]{catanese2011superficial}. This fact results in the following criterion proved in \cite[Proposition~5.7]{Wehler26Teich}.

\begin{prop}\label{prop: Tg exists}
	Suppose the complex space $\kc(\frg)$ is reduced. If $\kp_X$ is injective for every complex nilmanifold $X$ of type $(\frg, \Gamma)$ and the fibres of $\widetilde{\kp}$ are connected, then $\kp$ is injective and the Teichmüller space $\kt(\frg)$ admits a natural complex structure.
\end{prop}

We will use this criterion to establish a complex structure for all but two of the Lie algebras listed in the first item of Theorem \ref{thm: teichmüller and period maps dim 6}. For the remaining two cases, we will show directly that the map $\pi_X \colon \Kur(X) \to \kt(\frg)$ is a local homeomorphism and the Kuranishi spaces can be glued to a global complex structure on $\kt(\frg)$.

The main obstruction to this approach, and to the map $\pi_X$ being injective, is the fact that the dimension of the automorphism group may vary within the Kuranishi family of $X$. This is made explicit by Proposition \ref{thm: h^0 constant}, which we will prove next. Let $G$ be a simply connected nilpotent Lie group associated to $\frg$, and let $G_\IC$ be its universal complexification. Both of these groups act on the variety 
\[
\kv(\frg) = \{V \in \Gr(n, \frg_\IC) \, | \, [V, V] \subset V\} \subset \Gr(n, \frg_\IC),
\]
but a priori only $G$ preserves the subspace $\kc(\frg) \subset \kv(\frg)$.

\begin{lem}\label{lem: complexified orbits}
	Let $J \in \kc(\frg)$ be a complex structure. Then $G_\IC \cdot J = G \cdot J$.
\end{lem}
\begin{proof}
	Let $\frg^{1,0}$ be the subspace of $\frg_\IC$ corresponding to $J$ and let $\exp(\frg^{1,0}) = G^{1,0}$ denote the associated subgroup of $G_\IC$. Based on a result of Bia\l ynicki-Birula and Rosenlicht \cite{Rosenlicht62}, it was shown in \cite{Hasegawa26} that the multiplication map 
	\[
	G \times G^{1,0} \to G_\IC 
	\]
	is bijective. So any $g_\IC \in G_\IC$ can be written as $g_\IC= g g^{1,0}$ for $g \in G$ and $g^{1,0} \in G^{1,0}$. Hence, we have
	\[
	g_\IC \cdot \frg^{1,0} = g g^{1,0} \cdot \frg^{1,0} = g \cdot \Ad_{g^{1,0}}(\frg^{1,0}) = g \cdot \frg^{1,0}.
	\]
	Therefore, for every element in $\kc(\frg)$, the $G$-orbit coincides with the $G_\IC$-orbit.
\end{proof}

Following the notation of \cite{Meersseman19} and \cite{Wehler26Teich}, the function $h^0 \colon \kc(\frg) \to \IN$ maps $J$ to the dimension of the group of holomorphic right-translations on $(G,J)$. In other words, $h^0(J) = n - \dim_\IC G\cdot J$. For any lattice $\Gamma\subset G$, the value $h^0(J)$ also coincides with the dimension of the automorphism group of the complex nilmanifold $(\Gamma\backslash G, J)$.

\begin{prop}\label{prop: jumping criterion}
	Let $\frg$ be a real nilpotent Lie algebra admitting complex structures. Suppose $J_0 \in \kc(\frg)$ is a complex structure such that the function $h^0$ is not locally constant at $J_0$. Then the Teichmüller space $\kt(\frg)$ is not locally Hausdorff at $[J_0] \in \kt(\frg)$.
\end{prop}
\begin{proof}
	We will show that in every neighbourhood of $[J_0]$ in $\kt(\frg)$, there is a point $[J]\neq [J_0]$ such that $[J_0]$ and $[J]$ are inseparable. Throughout the proof we assume that $\kv(\frg)$ is reduced since passing to its reduction does not change the argument. Let $\kc_0$ be the connected component of $\kc(\frg)$ containing $J_0$, and let $\kv_0$ be an irreducible component of $\kv(\frg)$ containing $J_0$ such that $h^0$ is not locally constant at $J_0$ in $\kc_0 \cap \kv_0$. We denote by $R$ the image of the map
	\[
	\mu \colon G_\IC \times \kv_0 \to \kv_0 \times \kv_0
	\]
	sending $(g_\IC,V)$ to $(V, g_\IC\cdot V)$. Since $G_\IC$ is a complex unipotent group, the map $\mu$ is a morphism of algebraic varieties, and thus its image $R$ is a constructible subset of $\kv_0 \times \kv_0$ \cite[Section~I.8, Corollary~2]{Mumford99}, that is, a finite union of locally closed sets in the Zariski topology. Therefore the closure $\overline{R}$ of $R$ in the euclidean topology coincides with the Zariski closure of $R$ \cite[Section~I.10, Corollary~1]{Mumford99}, so the euclidean closure $\overline{R}$ is an irreducible subvariety of $\kv_0 \times \kv_0$. In particular, we obtain a morphism $p_1 \colon \overline{R} \to \kv_0$ by projecting to the first factor. Let $F_0$ be an irreducible component of the fibre $p_1^{-1}(J_0)$ containing the point $(J_0, J_0)$. Then we have
	\[
	\dim F_0 \geq \dim \overline R - \dim \kv_0 = k,
	\]
	where $k$ is the generic (maximal) orbit dimension of the $G_\IC$-action on $\kv_0$ \cite[Theorem~1.25]{Shafarevich13}. Since $\kc_0\cap \kv_0$ is an open subset of $\kv_0$, it follows from Lemma~\ref{lem: complexified orbits} that $k$ is also the maximal orbit dimension of the $G$-action on $\kc_0\cap \kv_0$. By our assumption on $h^0$ and our choice of $\kv_0$, the orbit $G \cdot J_0 = G_\IC \cdot J_0$ is a complex manifold of dimension strictly smaller than $k$. Let $U_0$ be the intersection of a euclidean neighbourhood of $(J_0,J_0) \in \kc_0 \times \kc_0$ with the smooth locus of $F_0$. Then $U_0$ is a complex manifold of dimension at least $k$, so there exist complex structures $J \in \kc_0$ arbitrarily close to $J_0$ such that $(J_0, J) \in \overline R \cap (\kc_0 \times \kc_0)$ but $J \not\in G \cdot J_0$.
	
	 It remains to show that for any such $J$, the classes $[J_0]$ and $[J]$ are inseparable in the connected component $\kc_0/G \subset \kt(\frg)$.  Since $(J_0, J) \in \overline R \cap (\kc_0 \times \kc_0)$, we can choose converging sequences $J_m\to J_0$ and $g_m\cdot J_m \to J$, with $J_m \in \kc_0$ and $g_m \in G$, again using Lemma \ref{lem: complexified orbits}. Hence, the sequence $[J_m]$ converges to $[J_0]$ and $[J]$ in $\kt(\frg)$, but since we showed that $[J_0]\neq [J]$, and we can choose $J$ arbitrarily close to $J_0$, this implies that the Teichmüller space is not locally Hausdorff at $[J_0]$.
\end{proof}

\subsection{The moduli space}\label{subsect: moduli criteria}

We will see that in most cases, at least up to dimension six, the function $h^0$ is constant. So in order to study the moduli space $\km(\frg, \Gamma)$ we need more refined criteria.

We begin by recalling the notion of an arithmetic group: a $\IQ$-algebraic group $\textbf{G}< \GL(m, \IC)$ is a subgroup defined by polynomials with rational coefficients in the coordinates of $\GL(m, \IC)$. If $R$ is a subring of $\IC$ we set $\textbf{G}(R) = \textbf{G} \cap \GL(m, R)$. A subgroup $\Lambda < \textbf{G}(\IQ)$ is said to be arithmetic if $\Lambda$ is commensurable with $\textbf{G}(\IZ)$, that is, the intersection $\Lambda \cap \textbf{G}(\IZ)$ is of finite index in both $\Lambda$ and $\textbf{G}(\IZ)$. 

Suppose $\frg$ has dimension $2n$. Every lattice $\Gamma$ induces a rational structure on $\frg$, i.e., a rational subalgebra $\frg_{\IQ}\subset \frg$ with $\frg_{\IQ} \otimes \IR = \frg$, by setting $\frg_\IQ = \langle \log \Gamma \rangle_\IQ$. With respect to a basis of $\frg_{\IQ}$ the automorphism group $\Aut(\frg_\IC)$ can be described by polynomials with rational coefficients, so $\textbf{G}= \Aut(\frg_\IC)$ is a $\IQ$-algebraic group in $\GL(2n, \IC)$ with $\textbf{G}(\IR) = \Aut(\frg)$ and $\textbf{G}(\IQ) = \Aut(\frg_\IQ)$. Throughout this section we will consider $\Aut(\frg)$ as a subgroup of $\GL(2n, \IR)$ with respect to a fixed basis of $\frg_{\IQ}$. Every element of $\Aut(\Gamma)$ extends uniquely to an automorphism of $G$, and by taking its differential we get an element of $\Aut(\frg)$. Therefore, we can consider $\Aut(\Gamma)$ as a discrete subgroup of $\textbf{G}(\IQ)$ and by \cite[Chapter~6, Theorem~6]{Segal83} the group $\Aut(\Gamma)$ is arithmetic.

Recall that another lattice $\Gamma' \subset G$ is commensurable with $\Gamma$ if and only if their associated rational structures coincide \cite[Theorem~5.1.12, p.~205]{Cor-Green}. Hence, for any lattice $\Gamma'$ commensurable with $\Gamma$ the group $\Aut(\Gamma') < \textbf{G}(\IQ)$ is also arithmetic with respect to the same $\IQ$-algebraic group.

\subsubsection{Properly discontinuous actions}
Our first goal will be to study when a complex structure on $\kt(\frg)$ descends to the quotient $\km(\frg, \Gamma)$.

The following classical result due to Cartan \cite{Cartan57} states that this hinges on the properness of the action by the lattice. More precisely, if $\Lambda$ is a group acting on a reduced complex space $Y$, then the complex structure on $Y$ descends to the quotient $Y/\Lambda$ if the $\Lambda$-action is properly discontinuous. In this case, the complex space $Y/\Lambda$ is normal if $Y$ is normal.

The group $\Aut(\Gamma)$ never acts properly discontinuously on $\kt(\frg)$ directly, since the inner automorphisms of $\Aut(\Gamma)$ always act trivially. So we will usually consider the group
\[
O_\Gamma = \im \left (\Aut(\Gamma) \to \Out(\frg)\right)
\]
and we can write $\km(\frg, \Gamma) = \kt(\frg)/O_\Gamma$. By abuse of notation we will sometimes write $\det \varphi$ for $\varphi \in \Out(\frg) = \Aut(\frg)/ \Inn(\frg)$ since the determinant of any inner automorphism is 1. We will denote by $\Out^1(\frg)$ the set of outer automorphisms whose determinant has absolute value $1$.
In many cases, the automorphism group $\Aut(\frg)$ acts transitively on $\kc(\frg)$ or unions of its connected components. In this case, $\kt(\frg)$ can be written as a union of connected components of the form $\Out(\frg)/H$ for some closed subgroup $H \subset \Out(\frg)$, and we have the following criterion.

\begin{prop}\label{prop: transitive action}
	Suppose $\kt(\frg)$ is a union of connected components
	\[
	\kt(\frg) = \bigsqcup_{i=1}^m \Out(\frg)/H_i,
	\]
	 where $H_i \subset \Out(\frg)$ is a closed subgroup for every $i=1,\dots,m$.
	\begin{enumerate}
		\item If $H_i$ is compact for every $i=1,\dots,m$, then $O_\Gamma$ acts properly discontinuously on $\kt(\frg)$.
		\item If $\Out^1(\frg)H_i = \Out(\frg)$ and $\Out^1(\frg) \cap H_i$ is compact for every $i=1,\dots,m$, then $O_\Gamma$ acts properly discontinuously on $\kt(\frg)$.
	\end{enumerate}
\end{prop}
\begin{proof}
	The first item follows immediately from the fact that any discrete subgroup of $\Out(\frg)$ acts properly discontinuously on the quotient $\Out(\frg)/H_i$ if $H_i$ is compact \cite[Lemma~3.1.1]{Wolf84}.
	
	For the second item, we first observe that for any automorphism $\varphi \in \Aut(\Gamma) \subset \Aut(\frg)$ we have $|\det \varphi| =1$. Indeed, since $\Aut(\Gamma)$ is arithmetic, the subgroup of integer matrices in $\Aut(\Gamma)$ is of finite index in $\Aut(\Gamma)$. Hence, $\det$ takes values in $\{\pm 1\}$ on a finite index subgroup of $\Aut(\Gamma)$, so $\im( |\det| \colon \Aut(\Gamma) \to \IQ_{>0})$ is a finite subgroup of $\IQ_{>0}$ and thus has to be $\{1\}$. 
	
	Hence, we have $O_\Gamma \subset \Out^1(\frg)$. The assumption $\Out^1(\frg)H_i = \Out(\frg)$ implies that the inclusion $\Out^1(\frg) \subset \Out(\frg)$ induces an $O_\Gamma$-equivariant homeomorphism
	\[
	\Out^1(\frg)/(\Out^1(\frg) \cap H_i) \to \Out(\frg)/H_i.
	\]
	So if $\Out^1(\frg) \cap H_i$ is compact, then by the same argument as in the first item, the group $O_\Gamma$ acts properly discontinuously on $\Out^1(\frg)/(\Out^1(\frg) \cap H_i) \cong \Out(\frg)/H_i$.
\end{proof}

Proposition \ref{prop: transitive action} will be used in Section \ref{sect: moduli dim six} to show that the moduli space for complex nilmanifolds of type $(\frh_{19}^-, \Gamma)$ and $(\frh_{26}^+, \Gamma)$ is a normal complex space.

\subsubsection{Non-Hausdorff criteria}
Next, we will provide two criteria to determine when the moduli space $\km(\frg, \Gamma)$ is non-Hausdorff, assuming throughout that the Teichmüller space $\kt(\frg)$ is Hausdorff.  

The proof of the following criterion goes along the same lines as the proof of \cite[Section~7.2, Theorem~2.4]{lange2012complex}.

\begin{prop}\label{prop: moduli non Hausdorff}
	 Let $\textbf{H}\subset \GL(2n, \IC)$ be a connected semisimple $\IQ$-algebraic group such that $H =\textbf{H}(\IR)$ is connected and contained in $\Aut(\frg)$. Suppose there exists a point $p_0 \in\kt(\frg)$ such that the identity component $H_0$ of the stabiliser for the $H$-action on $\kt(\frg)$ is non-compact. Let $N$ be the largest normal subgroup of $H$ contained in the identity component $H_0$.
	
	If $N \cap \Aut(\Gamma)$ is finite, or equivalently, $N \cap \GL(2n, \IZ)$ is finite, then $\km(\frg, \Gamma')$ contains non-closed points for any lattice $\Gamma' \subset G$ commensurable with $\Gamma$. In particular, $\km(\frg, \Gamma')$ is not Hausdorff.
\end{prop}
\begin{proof}
	We will show that there exists a point $q\in \kt(\frg)$ whose $\Aut(\Gamma)$-orbit in $\kt(\frg)$ is not closed. For such a point the class $[q] \in \km(\frg, \Gamma)$ is a non-closed point. Suppose the intersection $N \cap \Aut(\Gamma)$ is finite. Then the arithmetic subgroup $\Lambda =H \cap \Aut(\Gamma)$ of $H$ also intersects $N$ in finitely many points. Let $V$ be an open subset in $Y = H_{0}\backslash H$ with compact closure. Then the group $\Lambda$ acts on $Y$ from the right and there is a point $y \in Y$ such that $\Lambda \cdot y \cap V$ consists of infinitely many points \cite[Section~7.2, Lemma~2.2]{lange2012complex}.
	
	Let $\pi \colon H_0\backslash H \to \kt(\frg)$ be the induced orbit map at $p_0 \in \kt(\frg)$. The fibres of $\pi$ are discrete and thus meet the closure $\overline{V}$ only in finitely many points. Therefore, the set $\pi(\Lambda \cdot y \cap V)$ is also infinite. In particular, for $q= \pi(y)$ the intersection $\Aut(\Gamma)\cdot q \cap \pi(\overline V)$ is infinite, so it admits an accumulation point. Hence, the orbit $\Aut(\Gamma)\cdot q$ is not closed in $\kt(\frg)$ (cf. \cite[Section~7.2, Lemma~2.3]{lange2012complex}), and thus the point $[q] \in \km(\frg, \Gamma)$ is not closed.
\end{proof}

\begin{rem}\label{rem: tori and Kodaira manifolds}
	Proposition \ref{prop: moduli non Hausdorff} generalises the well-known fact that the moduli space of complex tori of dimension $n$ is not Hausdorff for $n\geq 2$ \cite[Corollary~3.3, p.~216]{lange2012complex}. Indeed, if $\frg \cong \IR^{2n}$ is the abelian Lie algebra of dimension $2n$, then $H = \SL(2n, \IR)$ satisfies the assumptions of Proposition \ref{prop: moduli non Hausdorff} (cf. \cite[Theorem~2.4, p.~214]{lange2012complex}).
	
	Similarly, applying Proposition \ref{prop: moduli non Hausdorff} for a Lie algebra underlying a Kodaira manifold of real dimension $2n \geq 6$ and with $H = \Sp(2n-2, \IR)$ is essentially a reformulation of the argument used in \cite{grantcharov2004deformations} to show that the moduli space of Kodaira manifolds is not Hausdorff.
	In the list of six-dimensional Lie algebras, this rules out the existence of a complex moduli space for complex nilmanifolds of type $(\frh_1, \Gamma)$ and $(\frh_3, \Gamma)$.
\end{rem}

We will see in Section \ref{sect: moduli dim six} that Proposition \ref{prop: moduli non Hausdorff} can also be applied to some but not all remaining Lie algebras. The Lie algebras not covered by either Proposition \ref{prop: transitive action} or Proposition \ref{prop: moduli non Hausdorff} will be studied via our last criterion.

Recall that the radical of an algebraic group is the identity component of its largest normal solvable subgroup.
Let us denote by $U$ the unipotent radical of the algebraic group $\Aut(\frg)$, that is, the set of unipotent elements in the radical of $\Aut(\frg)$. The unipotent radical $U$ is a connected unipotent normal subgroup of $\Aut(\frg)$ \cite[Section~11.21]{Borel91}. The group $\Lambda = \Aut(\Gamma) \cap U$ is an arithmetic subgroup of $U$, and since $U$ is unipotent $\Lambda$ is a co-compact lattice in $U$ \cite[Section~6.10]{Borel62}.

The following fact was proved by Shah in \cite[Proposition~3.2]{Shah91} for certain $\IQ$-algebraic groups. We recall a reformulation of his result presented as in \cite[Claim~3.11]{Entov23} for the unipotent radical $U\subset \Aut(\frg)$ and the lattice $\Lambda \subset U$: let $Q \subset U$ be a closed connected Lie subgroup, and let $S_Q \subset U$ be the identity component of the smallest $\IQ$-algebraic subgroup of $U$ containing $Q$. Then, for $u = \id_{\frg} Q \in U/Q$,  we have
\[
\overline{\Lambda \cdot u} = \Lambda S_Q \cdot u,
\]
where $\overline{\Lambda \cdot u}$ is the closure of the orbit $\Lambda \cdot u$ in $U/Q$.

\begin{prop}\label{prop: unipotent radical}
	Let $p_0 \in \kt(\frg)$ and let $A_0$ be the stabiliser of the $\Aut(\frg)$-action on $\kt(\frg)$ with respect to $p_0$. Let $\gotha_0$ and $\gothu$ denote the Lie algebra of $A_0$ and $U$, respectively.
	
	If $\gotha_0 \cap \gothu \subset \Der(\frg)$ is not an ideal, then the set of points in $\Aut(\frg) \cdot p_0 \subset \kt(\frg)$ whose $\Aut(\Gamma)$-orbit is not closed, is dense. Hence, the moduli space $\km(\frg, \Gamma)$ contains non-closed points and, in particular, is not Hausdorff. 
\end{prop}
\begin{proof}
	Consider the left-action of the unipotent radical $U \subset \Aut(\frg)$ on $\Aut(\frg)/A_0$, and let $U_{[\varphi]}$ be the identity component of the stabiliser of this action with respect to $[\varphi] = \varphi A_0 \in \Aut(\frg)/A_0$. Applying \cite[Claim~3.11]{Entov23} as stated above with $Q = U_{[\varphi]}$ we get
	\begin{equation}\label{eq: Shah orbit}
			\overline{\Lambda \cdot u} = \Lambda S_{U_{[\varphi]}} \cdot u \subset U/U_{[\varphi]}
	\end{equation}	
	 for $u = \id_{\frg} U_{[\varphi]}\in U/U_{[\varphi]}$. Now, the orbit map $U \to \Aut(\frg)/A_0$ with respect to $[\varphi] \in \Aut(\frg)/A_0$ descends to a $U$-equivariant continuous map 
	 \[
	 \pi_{\varphi} \colon U/U_{[\varphi]} \to \Aut(\frg)/A_0  = \Aut(\frg) \cdot p_0 \subset \kt(\frg).
	 \] 
	 Both $S_{U_{[\varphi]}}$ and $\Lambda$ act on $\kt(\frg)$ by restricting the $U$-action on $\kt(\frg)$, and for the element $p = \varphi \cdot p_0 \in \Aut(\frg) \cdot p_0$ the equality in \eqref{eq: Shah orbit} yields
	 \[
	  S_{U_{[\varphi]}} \cdot p =\pi_{\varphi}(S_{U_{[\varphi]}} \cdot u) \subset \pi_{\varphi}(\overline{\Lambda \cdot u}) \subset \overline{\pi_{\varphi}(\Lambda \cdot u)} = \overline{\Lambda \cdot p}.
	 \]
	 Therefore, we get
	 \[
	  S_{U_{[\varphi]}} \cdot p \subset \overline{\Lambda \cdot p} \subset \overline{\Aut(\Gamma) \cdot p} \subset \kt(\frg).
	 \]
	 Since $\Aut(\Gamma)$ is a countable group, this implies that the orbit $\Aut(\Gamma) \cdot p$ can only be closed if the algebraic group $S_{U_{[\varphi]}}$ acts trivially on $p$. The group $S_{U_{[\varphi]}}$ is the identity component of the smallest $\IQ$-algebraic subgroup of $U$ containing $U_{[\varphi]}$, so $S_{U_{[\varphi]}}$ acts trivially on $p$ if and only if $U_{[\varphi]} =S_{U_{[\varphi]}}$. We claim that if $\gotha_0 \cap \gothu$ is not an ideal of $\Der(\frg)$ then we can choose $\varphi$ such that $U_{[\varphi]} \subsetneq S_{U_{[\varphi]}}$.

	 Note that the Lie algebra of $U_{[\varphi]}$ is $\gothu_{[\varphi]} = \Ad_{\varphi}(\gotha_0 \cap \gothu)$, where $\Ad_{\varphi} \colon \Der(\frg) \to \Der(\frg)$ is the adjoint map at $\varphi \in \Aut(\frg)$. Suppose we have $\dim \gotha_0 \cap \gothu =k$, then taking $\varphi$ to $\gothu_{[\varphi]}$ defines a map
	 \[
	 \Aut(\frg) \to \Gr(k, \gothu),
	 \]
	whose image is positive-dimensional because we assumed that $\gotha_0 \cap \gothu$ is not an ideal. In particular, there is an element $\psi \in \Aut(\frg)$ such that $\gothu_{[\psi]}$ is not rational with respect to the rational structure $\gothu_{\IQ}$ on $\gothu$ induced by the lattice $\Lambda = \Aut(\Gamma) \cap U$. Therefore, we have $U_{[\psi]} \subsetneq S_{U_{[\psi]}}$ since the Lie algebra of $S_{U_{[\psi]}}$ is rational with respect to $\gothu_\IQ$. Hence, the closure of the $\Aut(\Gamma)$-orbit of $p' = \psi \cdot p_0$ contains the positive-dimensional orbit $S_{U_{[\psi]}} \cdot p'$, and thus $\Aut(\Gamma) \cdot p'$ is not closed.
	
	We have shown that the $\Aut(\Gamma)$-orbit of a point $p' =\psi \cdot p_0$ with $\psi \in \Aut(\frg)$ is not closed as soon as the subalgebra $\gothu_{[\psi]}$ is not rational. The set of automorphisms $\varphi$ for which $\gothu_{[\varphi]}$ is rational is a countable union of cosets of the form $\varphi N$, where $N$ is the stabiliser at $\gotha_0\cap\gothu$ for the $\Aut(\frg)$-action on $\Gr(k, \gothu)$. Since the complement of this union is dense in $\Aut(\frg)$, the set of points in $\Aut(\frg) \cdot p_0 \subset \kt(\frg)$ whose $\Aut(\Gamma)$-orbit is not closed is dense. This concludes the proof.
\end{proof}

\section{Teichmüller spaces for six-dimensional nilmanifolds} \label{sect: proof global period map dim 6}

In this section, we collect all the necessary results to prove Theorem \ref{thm: teichmüller and period maps dim 6}. As already mentioned, the period map $\kp$ was already studied for the Lie algebras $\frh_2, \frh_3, \frh_4, \frh_5, \frh_6, \frh_8, \frh_9$, and $\frh_{10}$ in \cite[Theorem~10.1]{Wehler26Teich}. After proving that the period map $\kp$ is injective for the Lie algebras $\frh_7$, $\frh_{11}$, $\frh_{12}$, $\frh_{19}^-$ and $\frh_{26}^+$, we will treat the remaining Lie algebras $\frh_{13}$, $\frh_{14}$, and $\frh_{16}$ in Sections \ref{subsect: h13 and h14} and \ref{sect: h16}.

\subsection{Kuranishi spaces and period maps}

Building on the classification of complex structures on six-dimensional Lie algebras given in \cite{ceballos2016invariant}, we first study the local structure of the Teichmüller spaces.

\begin{prop}\label{prop: Kuranishi spaces dim 6}
	Let $X=(M,J)$ be a complex nilmanifold of dimension three with underlying Lie algebra $\frg$. If $\frg \not\in \{ \frh_5, \frh_{15}, \frh_{16}\}$, then the following holds.
	
	\begin{enumerate}
	 \item The Kuranishi family of $X$ is smooth, universal, and the dimension of $\Kur(X)$ depends only on the underlying Lie algebra.
	\begin{center}
		\begin{tabular}{l c}
			\toprule
			\textup{Lie algebra} & $\dim \Kur(X)$\\
			\midrule
			$\frh_1$ & $9$\\
			$\frh_2,\frh_3,\frh_4$ & $4$\\
			$\frh_6,\frh_7$& $5$\\
			$\frh_8$ & $6$\\
			$\frh_9,\frh_{10},\frh_{11},\frh_{12}, \frh_{13},\frh_{14}$ & $3$\\
			$\frh_{19}^-,\frh_{26}^+$& $2$\\
			\bottomrule
		\end{tabular}
	\end{center}
	\item If $\frg \not\in \{\frh_{13}, \frh_{14}\}$, then the local period map $\kp_X \colon \Kur(X) \to \IP H^3(\frg, \IC)$ is an embedding.
	\end{enumerate}
\end{prop}
\begin{proof}
	
	By \cite[Theorem~4.1, Remark~4.2]{ceballos2016invariant} the Frölicher spectral sequence of $X$ degenerates at the first page for every complex structure on $\frh_1$, $\frh_3$, $\frh_6$, $\frh_7$, $\frh_8$, $\frh_9$, $\frh_{10}$, $\frh_{11}$, $\frh_{12}$ and $\frh_{19}^-$. Since complex nilmanifolds have trivial canonical bundle it follows from \cite[Theorem~3.3]{Anthes18} that in any of these cases the Kuranishi space of $X$ is smooth and the Kuranishi family is universal. Moreover, the local period map $\kp_X$ is an embedding in all of these cases by \cite[Lemma~5.6]{Wehler26Teich}. The same holds for every complex structure on the Lie algebras $\frh_2$ and $\frh_4$ by \cite[Theorem~10.1]{Wehler26Teich}.
	
	It remains to consider the Lie algebras $\frh_{13}$, $\frh_{14}$, and $\frh_{26}^+$.
	
	\begin{prooflist}
		\item[Case 1] $\frg \in \{\frh_{13}, \frh_{14}\}$. By \cite[Table~1]{ceballos2016invariant} every complex structure $J$ has a basis of $(1,0)$-forms $\omega^1,\omega^2,\omega^3\in\frg^{*1,0}$ with structure equations
		\[
		d\omega^1=0,\; d\omega^2=\omega^{1\bar1},\; d\omega^3=\omega^{12}+B\,\omega^{1\bar2}+c\,\omega^{2\bar1},
		\]
		where $B \in \IC$ and $c \in \IR\backslash\{0\}$. Let $z_1,z_2,z_3\in\frg^{1,0}$ be the dual basis. Then we compute
		\[
		\dbar z_1=\overline\omega^1\otimes z_2+B\overline\omega^2\otimes z_3,\;
		\dbar z_2=c\overline\omega^1\otimes z_3,\; \dbar z_3=0,
		\]
		and thus
		\[
		H^1(X,\Theta_X)=\langle[\overline\omega^1\otimes z_1+B\,\overline\omega^3\otimes z_3],[\overline\omega^2\otimes z_2-c\,\overline\omega^3\otimes z_3],[\overline\omega^2\otimes z_3]\rangle .
		\]
		 The non-zero Schouten brackets of the summands of the above representatives are
		\begin{align*}
			[\overline\omega^1\otimes z_1,\overline\omega^2\otimes z_2]&=\omega^{\overline{12}}\otimes[z_1,z_2]+\overline\omega^1\wedge(z_1\lrcorner d\overline\omega^2)\otimes z_2=-\omega^{\overline{12}}\otimes z_3,\\
			[\overline\omega^1\otimes z_1,\overline\omega^3\otimes z_3]&=\overline\omega^1\wedge(z_1\lrcorner d\overline\omega^3)\otimes z_3=-c\omega^{\overline{12}}\otimes z_3,\\
			[\overline\omega^2\otimes z_2,\overline\omega^3\otimes z_3]&=\overline\omega^2\wedge(z_2\lrcorner d\overline\omega^3)\otimes z_3=\overline B\omega^{\overline{12}}\otimes z_3 .
		\end{align*}
		Since $\omega^{\overline{12}}\otimes z_3=\dbar(\overline\omega^3\otimes z_3)$ is exact, every class in $H^1(X,\Theta_X)$ can be integrated to a one-parameter family of deformations, and thus $\Kur(X)$ is smooth of dimension 3. Moreover, the Kuranishi family is universal because $h^0$ is constant on $\Kur(X)$~\cite{wavrik1969obstructions}.
		
		\item[Case 2] $\frg = \frh_{26}^+$. According to \cite[Proposition~2.3]{Ugarte14Non-nilpotent} every complex structure $J$ has a basis of $(1,0)$-forms $\omega^1, \omega^2, \omega^3 \in \frg^{*1,0}$ with structure equations
		\[
		d\omega^1 =0, \; d\omega^2 = \omega^{13}+ \omega^{1\bar3}, \; d\omega^3 = i \omega^{1\bar 1} +\sigma i(\omega^{1\bar 2} - \omega^{2\bar 1}),
		\]
		where $\sigma \in \{1,-1\}$. If $z_1, z_2, z_3 \in \frg^{1,0}$ is the dual basis, then we have
		\[
		H^1(X,\Theta_X) = \langle [\overline\omega^3\otimes z_3], [\overline\omega^1\otimes z_2]\rangle.
		\]
		The only possibly non-zero Schouten-bracket is $[\overline\omega^3\otimes z_3, \overline\omega^1\otimes z_2]$, so we compute
		\[
		[\overline\omega^3\otimes z_3, \overline{\omega}^1\otimes z_2] = \overline\omega^1 \wedge(z_2 \lrcorner d\overline\omega^3) \otimes z_3 = i \sigma\overline\omega^1 \wedge \overline\omega^1 \otimes z_3 = 0.
		\]
		Therefore, $[\phi, \phi]=0$ for every $\phi \in \langle \overline\omega^3\otimes z_3,\, \overline{\omega}^1\otimes z_2\rangle$ and thus $\Kur(X)$ is smooth. As in the previous case, the Kuranishi family is also universal since $h^0$ is constant. In order to prove that $\kp_X$ is an embedding we have to show that every $d$-exact form in $\Wedge^{2,1}\frg^*$ is already $\dbar$-exact~\cite[Lemma~5.6]{Wehler26Teich}. Computing in the basis above we get
		\[
		\im d \cap \Wedge^{2,1}\frg^* = \langle \omega^{12\bar 1}, i\sigma\omega^{12\bar 2} - \omega^{13\bar 3} \rangle = \dbar(\Wedge^{2,0}\frg^*),
		\]
		so $\kp_X$ is an embedding.
	\end{prooflist}	
	Since $\Kur(X)$ is smooth in all of the above cases, its dimension coincides with the dimension of $H^{2,1}(\frg, J) \cong H^1(X,\Theta_X)$, for which we refer to an earlier version of \cite{ceballos2016invariant}, where these cohomology groups were explicitly computed for every pair $(\frg, J)$ (cf. \cite[Propositions~6.1-6.3]{ceballos2011classification}).
\end{proof}

	\begin{rem}\label{rem: dim C(g)}
	If $X= (M,J)$ is a complex nilmanifold with underlying Lie algebra $\frg$, then there is a small neighbourhood $U \subset \kc(\frg)$ of $J$ such that $U \cong \Kur(X) \times V$, where $V$ is an open neighbourhood around $0$ in a complement of $H^0(X, \Theta_X) \subset \frg^{1,0}$. In particular, $\kc(\frg)$ is smooth and of dimension $\Kur(X) + 3-h^0(J)$ for all Lie algebras listed in the first item of Proposition \ref{prop: Kuranishi spaces dim 6}. These dimensions coincide with the upper bound for $\dim_\IC \kc(\frg)$ given in \cite[Table~A.1]{salamon01}. For the remaining Lie algebras $\frh_5, \frh_{15}$ and $\frh_{16}$, the complex space is non-reduced by the criterion \cite[Proposition~3.4]{Wehler26Teich}, which causes $\frh_{15}$ and $\frh_{16}$ to have dimension one less than the upper bound in \cite[Table~A.1]{salamon01}.
\end{rem}

In order to apply Proposition \ref{prop: Tg exists} it remains to show that the fibres of the map $\widetilde{\kp} \colon \kc(\frg) \to \IP H^3(\frg, \IC)$ are connected.

\begin{prop}\label{prop: global period map dim 6}
	Let $\frg$ be a nilpotent Lie algebra of dimension six admitting complex structures. Then the fibres of $\widetilde \kp$ are connected. 
	
	Moreover, if $\frg \in \{\frh_{13}, \frh_{14}, \frh_{15}, \frh_{16}\}$, then $\dim \widetilde\kp^{-1}(\widetilde \kp(J)) > \dim G \cdot J$ for every $J \in \kc(\frg)$.
\end{prop}
\begin{proof}
	As already mentioned, it remains to treat the Lie algebras $\frh_7$, $\frh_{11}$, $\frh_{12}$, $\frh_{13}$, $\frh_{14}$, $\frh_{15}$, $\frh_{16}$, $\frh_{19}^-$, and $\frh_{26}^+$. Throughout this proof $\omega^1, \omega^2, \omega^3$ is always a basis of $(1,0)$-forms with respect to a complex structure $J$ on $\frg$, and with structure equations as in \cite[Table~1]{ceballos2016invariant}. We denote $\omega= \omega^1 \wedge \omega^2 \wedge \omega^3$, and we will describe the fibre 
	\[
	\widetilde\kp^{-1}(\widetilde \kp(J)) = \left\{ V \in \Gr(3, \frg_\IC^*)\, \middle| \, V \cap \overline V =0,\, \Wedge^3V \subset \langle \omega\rangle \oplus d\left(\Wedge^2\frg_\IC^*\right)\right\}
	\]
	for all of the above cases in terms of these $(1,0)$-forms. 
	\begin{prooflist}
		\item[Case 1.1]$\frg = \frh_7$. We have
		\[
		d\omega^1 =0, \; d\omega^2 = \omega^{1\bar1}, \; d\omega^3 = \omega^{12} + \omega^{1\bar2}.
		\]
		So we compute
		\[
		d\left(\Wedge^2\frg_\IC^*\right)= \langle \omega^{12\bar 1}- \omega^{1\overline{12}}, \omega^{12\bar 2}- \omega^{13\bar 1}, \omega^{1\overline{13}}- \omega^{2\overline{12}}, \omega^{12\bar 3}- \omega^{23\bar 1}- \omega^{3\overline{12}} + \omega^{1\overline{23}} \rangle.
		\]
		Since every complex structure on $\frg$ has a basis with the same structure equations as above, we only have to consider the exact forms $\omega^{12\bar 1}- \omega^{1\overline{12}}$ and $\omega^{12\bar 2}- \omega^{13\bar 1}$. Thus, we get
		\[
		\widetilde\kp^{-1}(\widetilde \kp(J)) = \left\{ \langle  \omega^1, \omega^2+ t\bar{\omega}^1, \omega^3 - t^2\overline{\omega}^1 + t\overline{\omega}^2 \rangle \in \Gr(3, \frg_\IC^*) \, | \, t \in \IC \right\},
		\]
		which is connected.
		
		\item[Case 1.2] $\frg \in \{ \frh_{11}, \frh_{12}\}$. The structure equations for both of these Lie algebras are of the form
		\[
		d\omega^1 =0,\, d\omega^2 = \omega^{1\bar 1}, \, d\omega^3 = \omega^{12} + B \omega^{1\bar 2} + |B-1|\omega^{2\bar 1},
		\]
		where both $B$ and $|B-1|$ are non-zero. So, for $c = |B-1|$ we have
		\[
		d\left(\Wedge^2\frg_\IC^*\right)= \langle\omega^{12\bar1}, \omega^{1\overline{12}},
		B\omega^{12\bar2}-\omega^{13\bar1}, (1-B)\omega^{12\bar2}+c\omega^{2\overline{12}},
		\omega^{1\overline{13}}-\overline B\omega^{2\overline{12}}, d\omega^{3\bar3} \rangle
		\]
		 Every complex structure on $\frg$ preserves the ascending central series \cite{rollenske09}, so every complex structure in the fibre $\widetilde\kp^{-1}(\widetilde \kp(J))$ admits a basis of $(1,0)$-forms $\xi^1, \xi^2, \xi^3$ such that
		\begin{equation*}
			\xi^1 = \omega^1,\;\xi^2 = \omega^2 + s \overline \omega^1,\; \xi^3 = \omega^3 + t_1 \overline \omega^1 + t_2 \overline \omega^2.
		\end{equation*}
		We claim that $\omega$ and $\xi = \xi^1 \wedge \xi^2 \wedge \xi^3$ are cohomologous in $H^3(\frg, \IC)$ if and only if $t_2 = s B$. So we compute
		\begin{align*}
			\xi &= \omega^1 \wedge (\omega^2 + s \overline \omega^1) \wedge (\omega^3 + t_1 \overline \omega^1 + t_2 \overline \omega^2)\\
			&= \omega - s \omega^{13\bar 1} + t_1\omega^{12\bar 1} + t_2 \omega^{12\bar 2} + st_2 \omega^{1\overline{12}}.
		\end{align*}
		 Therefore, $\xi = \xi^1 \wedge \xi^2 \wedge \xi^3$ is cohomologous to $\omega$ if and only if $t_2 = s B$ and thus the fibre 
		\[
		\widetilde\kp^{-1}(\widetilde \kp(J)) = \left\{ \langle \omega^1, \omega^2+ s \overline \omega^1, \omega^3 + t \overline \omega^1 + s B \overline \omega^2 \rangle \in \Gr(3, \frg_\IC^*) \, | \, s, t \in \IC \right\}
		\]
		is connected.
		\item[Case 1.3] $\frg \in \{ \frh_{19}^-, \frh_{26}^+\}$. The structure equations for a complex structure on these Lie algebras are
		\[
		d\omega^1 =0, \; d\omega^2 = \omega^{13} + \omega^{1\bar 3}, \; d\omega^3 = i \varepsilon \omega^{1\bar 1}+ \sigma i (\omega^{1\bar 2}- \omega^{2\bar 1}),
		\]
		where $\varepsilon \in \{0,1\}, \sigma \in \{1,-1\}$ and $\frg \cong \frh_{19}^-$ precisely if $\varepsilon =0$. We compute
		\begin{align*}
			d\left(\Wedge^2\frg_\IC^*\right) = \langle
			&\omega^{12\bar 1}, \omega^{1\overline{12}}, \omega^{13\bar 1}-\omega^{1\overline{13}},
			i\sigma\omega^{12\bar 2}-\omega^{13\bar 3}, i\sigma\omega^{2\overline{12}}+\omega^{3\overline{13}},\\
			&\omega^{13\bar 2}+\omega^{23\bar 1}-\omega^{1\overline{23}}-\omega^{2\overline{13}},\\
			&\sigma(\omega^{23\bar 1}-\omega^{13\bar 2}+\omega^{1\overline{23}}-\omega^{2\overline{13}})
			+\varepsilon(\omega^{13\bar 1}+\omega^{1\overline{13}})\rangle .
		\end{align*}
		Since every basis of $(1,0)$-forms for a complex structure on $\frg$ contains a closed form, we only have to consider the exact forms $\omega^{12\bar 1}$, $\omega^{1\overline{12}}$, $\omega^{13\bar 1}-\omega^{1\overline{13}}$, and $i\sigma\omega^{12\bar 2}-\omega^{13\bar 3}$. Therefore, the fibre $\widetilde\kp^{-1}(\widetilde \kp(J))$ is comprised of elements of the form
		\[
		\langle  \omega^1, \omega^2 -i\sigma st \overline \omega^1 -i \sigma t^2 \overline \omega^2 + t\overline \omega^3, \omega^3 +i(\sigma s + \varepsilon t) \overline \omega^1 + it\sigma \overline \omega^2 \rangle \in \Gr(3, \frg_\IC^*),
		\]
		where $s, t \in \IC$. Hence, $\widetilde\kp^{-1}(\widetilde \kp(J))$ is connected.
		
		\item[Case 2.1] $\frg \in \{\frh_{13}, \frh_{14}, \frh_{15}\}$. We claim that for every complex structure $J$ on $\frg$ with $h^0(J) =1$ the fibre $\widetilde\kp^{-1}(\widetilde \kp(J))$ is three-dimensional. This implies that $\kp$ is nowhere injective for $\frg \in \{\frh_{13},\frh_{14}\}$ and not generically injective for $\frg = \frh_{15}$ since the $G$-orbit $G\cdot J$ for such a complex structure is two-dimensional. We have the structure equations
		\[
		d\omega^1 =0,\; d\omega^2 = \omega^{1\bar 1}, \; d\omega^3 = \tau\omega^{12} + B \omega^{1\bar 2} + c\omega^{2\bar 1},
		\]
		with $\tau \in \{0,1\}, c\neq 0$ and $c \neq |B-\tau|$. It follows as in \emph{Case 1.2}, that every complex structure in the fibre $\widetilde\kp^{-1}(\widetilde \kp(J))$ can be described by a basis of the form
		\begin{equation*}
			\xi^1 = \omega^1, \; \xi^2 = \omega^2 + s \overline \omega^1, \;\xi^3 = \omega^3 + t_1 \overline \omega^1 + t_2 \overline \omega^2.
		\end{equation*}
		As in \emph{Case 1.2}, we compute
		\[
		\xi = \omega - s \omega^{13\bar 1} + t_1\omega^{12\bar 1} + t_2 \omega^{12\bar 2} + st_2 \omega^{1\overline{12}}.
		\]
		We have
		\[
		d\left(\Wedge^2\frg_\IC^*\right)= \langle \omega^{12\bar1}, \omega^{1\overline{12}}, \omega^{12\bar2}, \omega^{13\bar1}, \omega^{2\overline{12}}, \omega^{1\overline{13}}, d\omega^{3\bar3}\rangle,
		\]
		and thus $\xi$ is cohomologous to $\omega$ for every choice of $s, t_1$ and $t_2$. Similarly, one computes that $\dim \widetilde\kp^{-1}(\widetilde \kp(J)) > \dim G \cdot J$ if $\frg = \frh_{15}$ and $h^0(J) =2$.
		\item[Case 2.2] $\frg = \frh_{16}$. Since $h^0(J)=2$ for every complex structure on $\frh_{16}$, the orbit $G\cdot J$ has dimension one. We will show that for every choice of $s$, $t_1$ and $t_2$, the complex structure given by the basis 
		\begin{equation*}
			\xi^1 = \omega^1,\;\xi^2 = \omega^2 + s \overline \omega^1, \;\xi^3 = \omega^3 + t_1 \overline \omega^1 + t_2 \overline \omega^2
		\end{equation*}
		defines a complex structure in the fibre $\widetilde\kp^{-1}(\widetilde \kp(J))$. The structure equations for a complex structure $J$ on $\frh_{16}$ are
		\[
		d\omega^1 =0, \; d\omega^2 = \omega^{1\bar 1},\; d\omega^3 = \omega^{12} + B\omega^{1\bar 2}
		\]
		with $|B| =1 $ and $B \neq 1$. We compute $d(\omega^{2\bar 2} + \omega^{1\bar 3}) = (\overline B -1)\omega^{12\bar 1}$, $d(B\omega^{1\bar 3} + \omega^{2\bar 2}) = (1-B) \omega^{1\overline{12}}$, $d\omega^{3\bar 2}= \omega^{12\bar 2} - \omega^{13\bar 1}$, and $d\omega^{23} = -\omega^{13\bar 1}+ B\omega^{12\bar 2}$. Hence, the forms $\omega$ and $\xi$ are cohomologous and the fibres of $\kp$ always have positive dimension.\qedhere
	\end{prooflist}
\end{proof}

\subsection{The Lie algebras $\frh_{13}$ and $\frh_{14}$}\label{subsect: h13 and h14}

In this section, we will show that the Teichmüller spaces for $\frh_{13}$ and $\frh_{14}$ are complex manifolds, despite the fact that their period maps are far from being injective. We first recall two facts about complex structures on $\frg \in \{\frh_{13}, \frh_{14}\}$. By the classification proved in \cite{ceballos2016invariant}, we have $h^0\equiv 1$ for both Lie algebras, and moreover, the stabiliser of the $G$-action on $\kc(\frg)$ with respect to any $J \in \kc(\frg)$ is the centre $\kz G$. Furthermore, it was shown in \cite[Proposition~3.10]{rollenske09}, that every complex structure on $\frg$ preserves the ascending central series $0 \subset \kz \frg \subset \kz^2\frg \subset \frg$ on $\frg$, where $\dim \kz \frg =2$ and $\dim \kz^2\frg =4$.

\begin{lem}\label{lem: Tg hausdorff for h13 and h14}
	If $\frg = \frh_{13}$ or $\frg = \frh_{14}$, then the group $G/\kz G$ acts freely and properly on $\kc(\frg)$. In particular, the Teichmüller space $\kt(\frg)$ is Hausdorff.
\end{lem}
\begin{proof}
	Since the stabiliser at every point in $\kc(\frg)$ is $\kz G$ it remains to show that the action of $G/\kz G$ is proper. To this purpose, we choose a basis $e_1, \dots, e_6$ of $\frg$ compatible with the ascending central series. In particular, we have $\kz \frg = \langle e_5, e_6\rangle$ and $\frg/\kz \frg \cong \langle e_1, \dots, e_4 \rangle$. We consider $\kc(\frg)$ as a subspace of $\End(\frg) \cong \IR^{36}$ with respect to this basis. The action of $G/ \kz G$ on $\kc(\frg)$ is proper if for any converging sequences $J_n \to J \in \kc(\frg)$ and $g_n \cdot J_n \to J' \in \kc(\frg)$ with $g_n = \exp\left(\sum_{i=1}^4t_i(n)e_i\right)$, the sequence $(t_1(n), \dots, t_4(n))$ in $\IR^4$ has a convergent subsequence.
	
	 We will first show that $(t_1(n), t_2(n))$ has a convergent subsequence. Since every complex structure on $\frg$ preserves the centre $\kz\frg$, we obtain a map $\kc(\frg) \to \kc(\frh)$, where $\frh = \frg/\kz\frg$. The group $H = G/ \kz G$ acts on $\kc(\frh)$ and the sequences above descend to sequences $ \hat J_n \to \hat J$ and $\hat g_n \cdot \hat J_n \to  \hat J'$ in $\kc(\frh)$, where $\hat g_n = \exp( t_1(n)[e_1]+t_2(n)[e_2])$. The Lie algebra $\frh$ is the real Lie algebra underlying a Kodaira surface and it follows from the results of \cite{Borcea84}, that the group $H/\kz H = \exp\left(\langle [e_1], [e_2]\rangle\right)$ acts properly on $\kc(\frh)$. Therefore, the sequence $(t_1(n), t_2(n))$ has a convergent subsequence.
	
	Passing to a sequence $(t_1(n), \dots, t_4(n))$ where $(t_1(n), t_2(n))$ converges, we will show that $(t_3(n), t_4(n))$ has a convergent subsequence. Observe that
	\[
	g_n \cdot J_n = \exp(t_3(n)e_3 + t_4(n)e_4)\cdot (\exp(t_1(n)e_1 + t_2(n)e_2)\cdot J_n)
	\]
	because the image of $\ad_{e_3}$ and $\ad_{e_4}$ is in the centre of $\frg$. The sequence $I_n =\exp(t_1(n)e_1 + t_2(n)e_2)\cdot J_n$ has a convergent subsequence $I_{n_k} \to I$, and we have
	\[
	\exp(t_3(n_k)e_3 + t_4(n_k)e_4)\cdot I_{n_k} = I_{n_k}- t_3(n_k)[\ad_{e_3},I_{n_k}] - t_4(n_k)[\ad_{e_4}, I_{n_k}].
	\]
	Therefore, $t_3(n_k)[\ad_{e_3},I_{n_k}] +t_4(n_k)[\ad_{e_4}, I_{n_k}]$ has to converge as well. The two vectors $[\ad_{e_3},I]$ and $[\ad_{e_4},I]$ considered in $\End(\frg)$ are linearly independent, since no linear combination of $e_3$ and $e_4$ is a holomorphic vector field with respect to $I$. Thus, the sequence $(t_3(n_k), t_4(n_k))$ also converges.
\end{proof}

\begin{lem}\label{lem: Kur = Teich h13 and h14}
	Let $X$ be a complex nilmanifold of type $(\frg, \Gamma)$. If $\frg = \frh_{13}$ or $\frg = \frh_{14}$, then the map $\pi_X \colon \Kur(X) \to \kt(\frg)$ is injective.
\end{lem}
\begin{proof}
	According to Proposition \ref{prop: Kuranishi spaces dim 6}, $\Kur(X)$ is smooth and three-dimensional, so we can consider $\Kur(X)$ as an open ball around $0 \in H^1(X, \Theta_X) \cong \IC^3$.
	Seeking a contradiction, we assume that $\pi_X \colon \Kur(X) \to \kt(\frg)$ is not injective for any shrinking of $\Kur(X)$. Then there are sequences $(s_n)$ and $(t_n)$ converging to $0 \in \Kur(X)$ such that for large enough $n$, we have $s_n \neq t_n$ and the corresponding complex nilmanifolds $X_{s_n}$ and $X_{t_n}$ are biholomorphic via a right-translation $R_{g_n} \colon X_{s_n} \to X_{t_n}$. By Kuranishi's third theorem \cite{Kuranishi69} (cf. \cite[Theorem~5]{catanese2011superficial}), since $h^0$ is constant, there is a neighbourhood $U \subset \Diff(X)$ of the identity, such that for large enough $n$ no diffeomorphism in $U$ can be a biholomorphism from $X_{s_n}$ to $X_{t_n}$. We will derive a contradiction to this fact.
	
	Since every element of the centre $\kz G$ is an automorphism of $X_{t_n}$ and $X_{s_n}$ we can assume that $g_n = \exp(x_n)$, where $x_n \in \frg$ is in a complement $V$ of $\kz \frg \subset \frg$. Now, if $J_{s_n}$, $J_{t_n}$ and $J$ are the left-invariant complex structures on $X_{s_n}$, $X_{t_n}$ and $X =X_0$, respectively, then both $J_{s_n}$ and $g_n \cdot J_{s_n} = J_{t_n}$ converge to $J$ in $\kc(\frg)$. By Lemma~\ref{lem: Tg hausdorff for h13 and h14}, this implies that the sequence $(x_n)$ has a convergent subsequence $x_{n_k} \to x \in V$, and thus $\exp(x) \cdot J = J$. Since $G/\kz G$ acts freely on $\kc(\frg)$, this implies that $x =0$. Therefore, the sequence of right-translations $R_{g_{n_k}}$ converges to the identity of $X$. This contradicts Kuranishi's third theorem since these right-translations are also biholomorphisms between different fibres of the Kuranishi family.
	Hence, the map $\pi_X$ is injective.
\end{proof}

\begin{prop}\label{prop: Teichmüller h13 and h14}
	The Teichmüller spaces $\kt(\frh_{13})$ and $\kt(\frh_{14})$ are complex manifolds of dimension $3$.
\end{prop}
\begin{proof}
	As explained in Section \ref{sect: Teichmüller and moduli}, for $\frg \in \{\frh_{13}, \frh_{14}\}$ Lemma \ref{lem: Kur = Teich h13 and h14} allows us to glue the Kuranishi charts to a complex structure on $\kt(\frg)$, and since $\kt(\frg)$ is Hausdorff by Lemma \ref{lem: Tg hausdorff for h13 and h14}, it is a complex manifold. The dimension of $\kt(\frg)$ coincides with the dimension of the Kuranishi space of any nilmanifold of type $(\frg, \Gamma)$, so by Proposition~\ref{prop: Kuranishi spaces dim 6}, $\kt(\frg)$ is three-dimensional.
\end{proof}

\subsection{The Lie algebra $\frh_{16}$}\label{sect: h16}

Among the six-dimensional nilpotent Lie algebras admitting a complex structure, $\frh_{16}$ is the only Lie algebra where the function $h^0$ is constant and the complex space $\kc(\frh_{16})$ is singular. In fact, it is everywhere non-reduced because for every complex structure in $\kc(\frh_{16})$ there are $\dbar$-closed forms in $\Wedge^{2,1}\frh_{16}^*$, which are not $d$-closed \cite[Proposition~3.4]{Wehler26Teich}. Nevertheless, we have

\begin{prop}\label{prop: h16 Hausdorff}
	The Teichmüller space $\kt(\frh_{16})$ is Hausdorff.
\end{prop}
\begin{proof}
 Throughout the proof we will consider $\kc(\frh_{16})$ as a subspace of $\IP\left( \Wedge^3 (\frh_{16})_\IC^* \right)$ via the Plücker embedding \cite[Lemma~3.2]{Wehler26Teich}, that is, we identify
 \[
 \kc(\frh_{16}) = \left\{ [\omega] \in \IP\left( \Wedge^3 (\frh_{16})_\IC^*\right) \; \middle| \; d\omega =0, \,\omega \wedge \overline{\omega} \neq 0,\, \textup{$\omega$ decomposable} \right\}.
 \] 
  Let us denote by $H_{16}$ a simply connected Lie group associated to $\frh_{16}$. The topological space $\kt(\frh_{16}) = \kc(\frh_{16})/H_{16}$ is Hausdorff if and only if the set
	\[
	R = \left\{\left([\omega], h \cdot [\omega] \right) \in \kc(\frh_{16}) \times \kc(\frh_{16}) \; \middle| \; [\omega] \in \kc(\frh_{16}), h \in H_{16} \right\}
	\] 
	is closed. According to \cite{ceballos2016invariant}, every class in $\kc(\frh_{16})$ can be represented by a decomposable $3$-form $\omega = \omega^1 \wedge \omega^2 \wedge \omega^3$ such that 
	\[d\omega^1 =0, \, d \omega^2 = \omega^{1\bar 1}, \, d \omega^3 = \omega^{12} + B \omega^{1\bar 2}\]
	 with $B \neq 1$ and $|B| =1$. The $H_{16}$-orbit of such a form is
	 \[
	 H_{16} \cdot \omega = \left\{\omega(\lambda) = \omega + \lambda( \omega^{13\bar 1}- B \omega^{12\bar 2}) + \frac{\lambda^2}{2}(2B\omega^{1\overline{12}} - (B+1)\omega^{12\bar 1})\; \middle| \; \lambda \in \IC \right\}.
	 \]
	 More precisely, if $z_1,z_2, z_3$ is the dual basis of $\omega^1, \omega^2, \omega^3$, then we have $\omega(\lambda) = \exp(\lambda z_1 + \overline \lambda \overline z_1) \cdot \omega$. 
	 Note that for $\lambda \to \infty$ the class $[\omega(\lambda)] = [\frac{1}{\lambda^2}\omega(\lambda)] \in \IP\left( \Wedge^3 (\frh_{16})_\IC^* \right)$ converges to
	 \begin{equation}\label{eq: limit point}
	 	 [\tfrac{1}{2}(2B\omega^{1\overline{12}} - (B+1)\omega^{12\bar 1})] = [\omega^{1\bar 1} \wedge \tfrac{1}{2}(2B \overline \omega^2 + (B+1)\omega^2)],
	 \end{equation}
	 which is no longer an element of $\kc(\frh_{16})$ since the wedge product with its complex conjugate is zero. In order to show that $R$ is closed in $\kc(\frh_{16}) \times \kc(\frh_{16})$ we choose a sequence $([\omega_n], h_n \cdot [\omega_n])$ converging to $([\omega], [\omega'])$ as $n \to \infty$. Now, since every $H_{16}$-orbit in $\kc(\frh_{16})$ is one-dimensional, the orbits of the classes $[\omega_n]$ are still generated by translations through $\exp(\lambda z_1 + \overline \lambda \overline z_1)$ for sufficiently large $n$. So for $n$ large enough we can write $h_n \cdot [\omega_n] = \exp(\lambda_n z_1 + \overline \lambda_n \overline z_1) \cdot [\omega_n]$ for some $\lambda_n \in \IC$, and as above we denote this class by $[\omega_n(\lambda_n)]$. Moreover, the sequence $(\lambda_n)$ is bounded.
	 
	 Indeed, if the sequence $(\lambda_n)$ were to be unbounded, then the limit $[\omega'] \in \kc(\frh_{16})$ of the sequence $([\omega_n(\lambda_n)])$ would be the class in \eqref{eq: limit point}, which is no longer in $\kc(\frh_{16})$. Hence, we can choose a subsequence of $(\lambda_n)$ converging to $\mu \in \IC$, and thus the sequence  $([\omega_n], h_n \cdot [\omega_n])$ converges to $([\omega], \exp(\mu z_1+ \overline \mu \overline z_1) \cdot [\omega]) \in R$. Therefore, the set $R$ is closed and $\kt(\frh_{16})$ is Hausdorff.
\end{proof}

\subsection{Proof of Theorem \ref{thm: teichmüller and period maps dim 6}}
The claims about the complex space $\kc(\frg)$ were already explained in Remark \ref{rem: dim C(g)}.
The results about $\kt(\frg)$ follow for all Lie algebras not isomorphic to $\frh_5, \frh_{13}, \frh_{14}, \frh_{15}$, or $\frh_{16}$ by applying Proposition \ref{prop: Tg exists}. The assumptions of Proposition \ref{prop: Tg exists} were verified for these Lie algebras in the Propositions~\ref{prop: Kuranishi spaces dim 6} and~\ref{prop: global period map dim 6}. The table of dimensions for $\kt(\frg)$ coincides with the one given in Proposition \ref{prop: Kuranishi spaces dim 6} since in all cases the Teichmüller space is locally isomorphic to the Kuranishi space.

For $\frg \in \{ \frh_5, \frh_{15}\}$, we observe from the classification in \cite{ceballos2016invariant} that the function $h^0$ is not constant on $\kc(\frg)$ and thus $\kt(\frg)$ is not Hausdorff by Proposition~\ref{prop: jumping criterion}. Moreover, by Proposition~\ref{prop: global period map dim 6} the period map is not injective for $\frh_{13}$, $\frh_{14}$, $\frh_{15}$, and $\frh_{16}$, while the non-injectivity for $\frh_5$ was proved in \cite[Example~6.8]{Wehler26Teich}. 

The Teichmüller spaces $\kt(\frh_{13})$ and $\kt(\frh_{14})$ are complex manifolds of dimension three by Proposition \ref{prop: Teichmüller h13 and h14}.
Finally, $\kt(\frh_{16})$ is Hausdorff by Proposition~\ref{prop: h16 Hausdorff}. This concludes the proof of Theorem \ref{thm: teichmüller and period maps dim 6}.\qed

\section{Moduli spaces for six-dimensional nilmanifolds}\label{sect: moduli dim six}

In this section, we will prove Theorem \ref{thm: coarse moduli space dim 6}. Before we get to the individual cases, we first observe that, apart from tori, nilmanifolds of type $(\frh_7, \Gamma)$ are the only ones where the subsets of biholomorphic complex nilmanifolds in the Teichmüller space are possibly dense.

\begin{prop}\label{prop: dense orbits}
	Let $\frg$ be a six-dimensional nilpotent Lie algebra admitting complex structures. If there exists a lattice $\Gamma$ and a point $p \in \kt(\frg)$ such that the orbit $\Aut(\Gamma) \cdot p \subset \kt(\frg)$ is dense, then $\frg$ is abelian or $\frg = \frh_7$.
\end{prop}
\begin{proof}
	It was shown in \cite{rollenske09} that every six-dimensional non-abelian Lie algebra $\frg$, except for $\frh_7$, admits what is called a stable torus bundle series. This means that every complex nilmanifold $X$ of type $(\frg, \Gamma)$ with $\frg \neq \frh_7$ has the structure of an iterated holomorphic torus bundle. Furthermore, it follows from the classification given in \cite[Theorem~B]{rollenske09} that each of these bundles has an elliptic curve as its base or fibre. Therefore, we can define a continuous map 
	\[
	 \varphi \colon \kt(\frg) \to \IC,
	\]
	where we consider $\IC$ as the moduli space of elliptic curves. Moreover, this map is $\Aut(\Gamma)$-periodic and not constant. Thus, no $\Aut(\Gamma)$-orbit in $\kt(\frg)$ can be dense since each of them is contained in a fibre of $\varphi$.
\end{proof}

 The results for the moduli spaces of nilmanifolds of type $(\frh_1, \Gamma)$ and $(\frh_3, \Gamma)$ were already discussed in Remark \ref{rem: tori and Kodaira manifolds}. The following table provides an overview which criteria will be used to study the remaining moduli spaces.
\begin{center}
	\begin{tabular}{l c}
		\toprule
		\textup{Lie algebra} & \textup{Proposition} \\
		\midrule
		$\frh_2, \frh_4$ &  \ref{prop: moduli non Hausdorff}\\
		$ \frh_k, \; 5 \leq k \leq 16$ &  \ref{prop: unipotent radical}\\
		$\frh_{19}^-, \frh_{26}^+$ & \ref{prop: transitive action}\\
		\bottomrule\\
	\end{tabular}
\end{center}

Verifying the assumptions of these Propositions in each case comes down to a very explicit computation in terms of a basis for $\frg$. In more geometric terms, Propositions~\ref{prop: moduli non Hausdorff} and \ref{prop: unipotent radical} correspond to the following phenomena causing the moduli spaces to be non-Hausdorff.
\begin{itemize}
\item Complex nilmanifolds of type $(\frh_2, \Gamma)$ and $(\frh_4, \Gamma)$ are bundles of elliptic curves over a complex two-dimensional torus, and a generic two-dimensional torus occurs as the base of such a bundle. So the induced map to the moduli space of two-dimensional complex tori, which has a dense subset of non-closed points, is dominant.

\item Every complex nilmanifold of type $(\frg, \Gamma)$ with $\frg = \frh_k$ for $5 \leq k \leq 16$ has the structure of an iterated holomorphic torus bundle with non-trivial twist deformations, that is, deformations preserving both the fibres and base of this bundle (cf. \cite{Hoefer} or \cite[Section~6]{Wehler26Teich}). Such deformations do not exist for complex nilmanifolds of type $(\frh_2, \Gamma)$ and $(\frh_4, \Gamma)$.

\end{itemize}

Throughout this section $e_1, \dots, e_6$ will always denote a basis of $\frg$ with dual basis $e^1,\dots, e^6$ and any subgroup of $\Aut(\frg)$ or $\Der(\frg)$ will be described as a matrix group with respect to this basis. For indecomposable Lie algebras of dimension six the group $\Aut(\frg)$ was described in \cite{Magnin07}, while the decomposable case is covered by \cite{Fisher13}, where the authors describe $\Aut(\frg)$ for all Lie algebras up to dimension five.

\subsection{The Lie algebras $\frh_2$ and $\frh_4$} 

In order to apply Proposition \ref{prop: moduli non Hausdorff} to prove that $\km(\frg, \Gamma)$ contains non-closed points we have to show the following: there is a group $H \subset \Aut(\frg)$ making up the real points of a $\IQ$-algebraic group with respect to a basis for $\langle\log(\Gamma)\rangle_\IQ$, and a point $p_0 \in \kt(\frg)$ such that the identity component $H_0$ of the stabiliser at $p_0$ is non-compact, and the largest normal subgroup $N$  of $H$ contained in $H_0$ intersects $\GL(6, \IZ)$ in finitely many points.

\subsubsection{The Lie algebra $\frh_2$} The non-zero brackets of $\frh_2$ are
\[
[e_1, e_2] = e_5, \qquad [e_3, e_4] = e_6,
\]
that is, $\frh_2$ is the direct sum of two 3-dimensional Heisenberg Lie algebras. It was proved in \cite{deGraaf07}, that the Lie algebra $\frh_2$ admits several rational structures. We will first show that Proposition \ref{prop: moduli non Hausdorff} can be applied to lattices $\Gamma$ whose associated rational structure $\langle \log(\Gamma)\rangle_\IQ$ is isomorphic to $\langle e_1, \dots, e_6 \rangle_\IQ$. The group
\[
H = \left\{\begin{pmatrix}
	A  & 0 & 0 &0\\
	0 & B & 0 & 0\\
	0 & 0 & 1 & 0\\
	0& 0 & 0&  1
\end{pmatrix}
\middle | A, B \in \SL(2, \IR)
\right\}
\]
is a $\IQ$-algebraic subgroup of $\Aut(\frh_2)$ with respect to the basis $e_1, \dots,e_6$. Let $J_0$ be the complex structure on $\frh_2$ corresponding to the form
\[
\omega_0 = (e^1+ (1+i)e^3) \wedge (e^2+ (1+i)e^4) \wedge (-\tfrac{1}{2}e^5+ie^6).
\]
The claim that $\km(\frh_2, \Gamma)$ contains non-closed points now follows directly from Proposition \ref{prop: moduli non Hausdorff} and the following lemma.

\begin{lem}\label{lem: stabiliser h2}
	The stabiliser $H_0$ is the subgroup
	\[
	H_{0} = \left\{\begin{pmatrix}
		A  & 0 & 0 &0\\
		0 & A & 0 & 0\\
		0 & 0 & 1 & 0\\
		0& 0 & 0&  1
	\end{pmatrix}
	\middle | A \in \SL(2, \IR)
	\right\} \subset H.
	\]
	In particular, $H_{0}$ is non-compact and the intersection of $\GL(6, \IZ)$ with the largest normal subgroup of $H$ contained in $H_{0}$ is finite.
\end{lem}
\begin{proof}
	The structure equations of the basis $e^1, \dots, e^6$ are $de^5 = -e^{12}$ and $de^6 = -e^{34}$, so we get
	\[
	\im \left( d \colon \Wedge^2 \frh_2^* \to \Wedge^3\frh_2^*\right) = \langle e^{123}, e^{124}, e^{134}, e^{234}, e^{126}-e^{345} \rangle.
	\]
	In particular, a matrix $C \in H$ acts trivially on the class $[\omega_0] \in \IP H^3(\frh_2,\IC)$ if and only if $C^*\omega_0 = \lambda \omega_0$ for some $\lambda \in \IC$. If 
	\[
	C=
	\begin{pmatrix}
		A  & 0 & 0 &0\\
		0 & B & 0 & 0\\
		0 & 0 & 1 & 0\\
		0& 0 & 0&  1
	\end{pmatrix},
	\]
	this is equivalent to the $4\times 4$ matrix $\begin{pmatrix}
		A & 0\\
		0 & B
	\end{pmatrix}$ being in the image of the embedding $\GL(2, \IC) \hookrightarrow \GL(4, \IR)$, which is equivalent to $A=B$. Let $N$ be a normal subgroup of $H \cong \SL(2, \IR) \times \SL(2, \IR)$. Since the diagonal of $\SL(2, \IR) \times \SL(2, \IR)$ is not normal and the only non-trivial normal subgroup of $\SL(2, \IR)$ is $\{I_2, -I_2\}$, we get that the largest normal subgroup of $H$ contained in $H_0$ is finite.
\end{proof}

According to the classification proved in \cite{deGraaf07} the rational structures on $\frh_2$ are parametrised by $\IQ_{>0}/(\IQ^*)^2$.
Concretely, every rational structure on $\frh_2$ can be described, up to isomorphism, by a basis $
f_1= e_1+ e_3, \, f_2 = \sqrt{d}(e_1-e_3), \, f_3 = e_2+e_4,\, f_4= \sqrt{d}(e_2-e_4), \, f_5= e_5+e_6, \, f_6 = \sqrt{d}(e_5-e_6)$, where $d \in \IZ$ is a positive square free integer. For $d=1$, this basis defines the same rational subalgebra as $e_1, \dots, e_6$. 
So, for every lattice $\Gamma$ there is a square free integer $d$ such that $\langle\log(\Gamma)\rangle_\IQ \cong \langle f_1, \dots, f_6\rangle_\IQ$. Let $P_d$ be the base change matrix from $e_1,\dots, e_6$ to $f_1, \dots, f_6$. Then, we can apply Proposition \ref{prop: moduli non Hausdorff} to the group $P_d H P_d^{-1}$, and as in Lemma \ref{lem: stabiliser h2} the stabiliser $P_d H_0 P_d^{-1}$ is non-compact and intersects $N$ in finitely many points.
This proves that $\km(\frh_2, \Gamma)$ contains non-closed points for any lattice.

The following example shows that the number of connected components of the moduli space $\km(\frh_2, \Gamma)$ varies with the lattice (see Remark \ref{rem: orientation}).

\begin{exam}\label{exam: lattice dependence}
	 An associated simply connected Lie group $H_2$ to the Lie algebra $\frh_2$ can be described as a product $H_2= H \times H$, where
	\[
	H = \left\{
	\begin{pmatrix}
		1 & x_1 & x_3\\
		0 & 1 & x_2\\
		0 & 0 & 1
	\end{pmatrix}
	\middle | x_i \in \IR
	\right\}
	\]
	is the three-dimensional Heisenberg group. Both $\Lambda = H \cap \GL(3, \IZ)$ and
	\[
	\Lambda' = \left\{
	\begin{pmatrix}
		1 & a_1 & \tfrac{1}{2}a_3\\
		0 & 1 & a_2\\
		0 & 0 & 1
	\end{pmatrix}
	\middle | a_i \in \IZ
	\right\}
	\]
	are lattices in $H$, and thus $\Gamma = \Lambda \times \Lambda$ and $\Gamma' = \Lambda \times \Lambda'$ are lattices in $H_2$. We fix an orientation on $M= \Gamma\backslash H_2$ and $M' = \Gamma' \backslash H_2$ by fixing a basis of $\frh_2 = \frh \oplus \frh$ in $\log(\Gamma)$. Swapping the two factors of $\Gamma = \Lambda \times \Lambda$ induces an orientation-reversing diffeomorphism of $M$. We will show that no such diffeomorphism exists on $M'$. Let $f$ be a diffeomorphism of $M'$. Then $f$ induces an automorphism of $\Gamma'$, which uniquely extends to an automorphism $\varphi$ of $H_2$ \cite{Vinberg00}, and $f$ reverses the orientation of $M'$ precisely if $d\varphi \colon \frh_2 \to \frh_2$ has negative determinant. Every Lie algebra automorphism of $\frh_2 = \frh\oplus \frh$ either preserves or swaps the subalgebras $\frh \oplus \kz \frh$ and $\kz\frh \oplus \frh$. If $d\varphi$ preserves $\frh \oplus \kz \frh$, then its determinant has to be positive, so $f$ preserves the orientation. If $d\varphi$ interchanges the factors, then $\varphi$ induces an isomorphism between $\Lambda \times \IZ$ and $\Lambda' \times \IZ$, but the abelianisation of $\Lambda'$ has torsion and the one of $\Lambda$ does not. Therefore, $M'$ does not admit an orientation-reversing diffeomorphism.
\end{exam}

\subsubsection{The Lie algebra $\frh_4$} The non-zero brackets of $\frh_4$ are
\[
[e_1, e_2] = e_5, \qquad [e_1, e_4] = [e_2, e_3] = e_6. 
\] 
In the notation of \cite{Magnin07} this is the Lie algebra $M4$. Again, from the classification proved in \cite{deGraaf07} it follows that the Lie algebra $\frh_4$ has a unique rational structure, so it suffices to work in the basis $e_1, \dots, e_6$. To show that the moduli space $\km(\frh_4, \Gamma)$ is not Hausdorff, one uses the same argument as for the Lie algebra $\frh_2$ applied to the group
\[
H = \left\{\begin{pmatrix}
	a_{11}  & a_{12} &  & &  &\\
	a_{21} & a_{22} &  &  &\textup{\huge 0} &\\
	&  & a_{11} & -a_{12} &&\\
	&  & -a_{21}& a_{22} &&\\
	&\textup{\huge 0}&&&1&0\\
	&&&&0&1
\end{pmatrix}
\middle | \begin{pmatrix}
	a_{11} & a_{12}\\
	a_{21} & a_{22}
\end{pmatrix} \in \SL(2, \IR)
\right\},
\]
and the complex structure corresponding to 
\[
\omega_0 = (e^1-ie^3+e^4)\wedge (e^2+ie^4)\wedge(2e^5+ie^6).
\]

\subsection{The Lie algebras $\frh_5, \frh_6, \frh_7, \frh_8, \frh_9, \frh_{10}, \frh_{11}, \frh_{12}, \frh_{13}, \frh_{14}, \frh_{15}$, and $\frh_{16}$}

We will use Proposition~\ref{prop: unipotent radical} to show that $\km(\frg, \Gamma)$ has non-closed points for all of these Lie algebras. Recall that $\gothu$, which is the largest ideal of $\Der(\frg)$ consisting of nilpotent elements~\cite{Chevalley55}, denotes the Lie algebra of the unipotent radical of $\Aut(\frg)$, and $\gotha_0$ denotes the Lie algebra of the stabiliser of the $\Aut(\frg)$-action with respect to a point $p_0 \in \kt(\frg)$. More precisely, if $p_0 = [J_0]$, then 
\[
\gotha_0 = \{ \delta \in \Der(\frg) \, | \, \delta \circ J_0 = J_0 \circ \delta\} + \{\ad_X \, | \, X \in \frg\}.
\]
We will also denote by $I_{ij}$ the $6\times6$ matrix with the single non-zero entry $1$ at $(i,j)$.

Concretely, we will show that there is a complex structure $J_0$ on $\frg$ such that $\gotha_0 \cap \gothu$ is not an ideal in $\Der(\frg)$ for the Lie algebra $\gotha_0$ of the stabiliser at $p_0 = [J_0] \in \kt(\frg)$.
For $\frg \in \{\frh_6,\frh_7, \frh_8, \frh_9, \frh_{10}\}$ this implies that the set of non-closed points in $\km(\frg, \Gamma)$ is dense because $\Aut(\frg)$ acts transitively on $\kt(\frg)$ in these cases \cite{ceballos2016invariant}. 

\subsubsection{The Lie algebra $\frh_5$}
The non-zero brackets of $\frh_5$ are
\[
[e_1, e_3]= [e_4,e_2] = e_5, \qquad [e_1, e_4] = [e_2,e_3] = e_6.
\]
Consider the complex structure given by $J_0(e_1) = e_2, J_0(e_3) = -e_4, J_0(e_5) = e_6$. Then the stabiliser $\gotha_0 \subset \Der(\frh_5)$ with respect to $p_0 =[J_0] \in \kt(\frh_5)$ consists of matrices of the form
\[
\begin{pmatrix}
	a_{11} & a_{12} & 0 &0& 0 & 0 \\
	-a_{12} & a_{11} & 0& 0& 0 & 0 \\
	0 & 0 & a_{33} & a_{34} &0 & 0 \\
	0 & 0 & -a_{34} & a_{33} &0 & 0 \\
	a_{51} & a_{52} & a_{53} & a_{54} & a_{11}+a_{33} & a_{12}+a_{34} \\
	-a_{52} & a_{51} & a_{63} & a_{64} & -a_{12}-a_{34} &a_{11}+a_{33}
\end{pmatrix}.
\]
The element $I_{53}$ is contained in $\gotha_0 \cap \gothu$ and $I_{31}+I_{42}$ is a derivation of $\frh_5$.  Computing their bracket in $\Der(\frh_5)$ we get
\[
[I_{53}, I_{31}+I_{42}] = I_{51}\not\in \gotha_0 \cap \gothu.
\]
Hence, $\gotha_0 \cap \gothu \subset \Der(\frh_5)$ is not an ideal.
\subsubsection{The Lie algebra $\frh_6$}
The non-zero brackets of $\frh_6$ are
\[
[e_1, e_2] = e_5, \qquad [e_1, e_3] = e_6.
\]
For the complex structure $J_0$ given by $J_0(e_1)= e_4, J_0(e_2) =e_3, J_0(e_5) = e_6$, the stabiliser $\gotha_0 \subset \Der(\frh_6)$ with respect to $p_0 =[J_0] \in \kt(\frh_6)$ consists of matrices of the form
\[
 \begin{pmatrix}
	a_{11} & 0 & 0 & 0 & 0&0 \\
	0 & a_{22} & a_{23} & 0 & 0 & 0 \\
	0 & -a_{23} & a_{22}& 0 & 0 & 0 \\
	0 & 0 & 0 & a_{11}& 0 & 0 \\
	a_{51} & a_{63} & a_{53} & a_{54} & a_{11}+a_{22} & a_{23} \\
	a_{61} & -a_{53} & a_{63} & a_{64} & -a_{23}& a_{11}+ a_{22}
\end{pmatrix}.
\]
We have $I_{61}-I_{54}\in \gotha_0 \cap \gothu$, and $I_{42}$ is a derivation of $\frh_6$. Computing their bracket, we get
\[
[I_{61}-I_{54}, I_{42}] = -I_{52} \not\in \gotha_0 \cap \gothu.
\]
Hence, $\gotha_0 \cap \gothu$ is not an ideal.

\subsubsection{The Lie algebra $\frh_7$}\label{subsubsect: h7}
The non-zero brackets of $\frh_7$ are
\[
[e_1, e_2] = e_4, \qquad [e_1, e_3] = e_5, \qquad [e_2, e_3] = e_6.
\]
 The stabiliser $\gotha_0 \subset \Der(\frh_7)$ for the complex structure $J_0(e_1)=e_2, J_0(e_3) = e_4, J_0(e_5)=e_6$ consists of matrices of the form
\[
\begin{pmatrix}
	a_{11} & a_{12} & 0& 0 & 0 & 0 \\
	-a_{12} & a_{11} &0 & 0 & 0 & 0 \\
	a_{31} & a_{32} &2a_{11} & 0 & 0 & 0 \\
	a_{41} & a_{42} &0 & 2a_{11} & 0 & 0 \\
	a_{51} & a_{52} & a_{42} -2a_{31} & a_{32} & 3a_{11}&a_{12} \\
	-a_{52} & a_{51} & -a_{41}-2a_{32} & -a_{31} & -a_{12} & 3a_{11}
\end{pmatrix}.
\]
Thus, we have $I_{52}-I_{61} \in \gotha_0 \cap \gothu$, and for $I_{12}+I_{56} \in \Der(\frh_7)$ we get
\[
[I_{52}-I_{61}, I_{12}+I_{56}] = -I_{62}+I_{51} \not\in \gotha_0 \cap \gothu.
\]

\subsubsection{The Lie algebra $\frh_8$}
The non-zero brackets of $\frh_8$ are
\[
 [e_1, e_2] = e_6.
\]
Again, we define $J_0(e_1)=e_2, J_0(e_3) = e_4, J_0(e_5)=e_6$ and $\gotha_0$ is made up of matrices of the form
\[
\begin{pmatrix}
	a_{11} & a_{12} & 0 & 0 &0 & 0 \\
	-a_{12} & a_{11} & 0 & 0 &0&0 \\
	a_{31}& a_{32} & a_{33} & a_{34} & 0 & 0 \\
	-a_{32}& a_{31} & -a_{34} & a_{33} & 0 & 0 \\
	a_{51} &a_{52} & a_{53} & a_{54} & 2a_{11} & 0 \\
	a_{61} & a_{62} & -a_{54} & a_{53} & 0 & 2a_{11}
\end{pmatrix}.
\]
For $I_{51} \in \gotha_0 \cap \gothu$ and $I_{35} \in \Der(\frh_8)$, we have
\[
[I_{51}, I_{35}] = -I_{31} \not\in \gotha_0 \cap \gothu.
\]
\subsubsection{The Lie algebra $\frh_9$}
The non-zero brackets of $\frh_9$ are
\[
[e_1, e_2] = e_5, \qquad [e_1, e_4]= [e_2, e_5] =e_6.
\]
We define $J_0(e_1)= e_2, J_0(e_3)= e_6, J_0(e_4) = e_5$. The stabiliser $\gotha_0$ consists of matrices of the form
\[
\begin{pmatrix}
	a_{11} & 0 & 0 & 0 & 0&0 \\
	0 & a_{11} & 0 &0&0 & 0 \\
	a_{31} & a_{32} & 3a_{11} & a_{34} & 0 & 0 \\
	a_{41} & \tfrac12 a_{34} & 0 & 2a_{11} & 0 & 0 \\
	\tfrac12 a_{34} -a_{65} & a_{52} & 0 & 0 & 2a_{11} & 0 \\
	a_{61} & a_{62} & 0 & a_{52} -a_{41} & a_{65} & 3a_{11}
\end{pmatrix}.
\]
Therefore, we have $I_{41}-I_{64}\in \gotha_0 \cap \gothu$, and for $I_{12}-I_{54} \in \Der(\frh_9)$ we get
\[
[I_{41}-I_{64}, I_{12}-I_{54}] = I_{51}+I_{42} \not\in \gotha_0 \cap \gothu.
\]
\subsubsection{The Lie algebra $\frh_{10}$}
The non-zero brackets of $\frh_{10}$ are
\[
[e_1, e_2] = e_4, \qquad [e_1, e_3]= e_5, \qquad [e_1,e_4] = e_6.
\]
For $J_0(e_1) =e_2, J_0(e_3) =e_4, J_0(e_5) = e_6$, the stabiliser is comprised of the matrices
\[
\begin{pmatrix}
	a_{11} & 0 & 0&0 & 0 & 0 \\
	0 & a_{11} & 0&0 &0 &0 \\
	a_{31} & a_{32} &2a_{11} & 0 & 0 & 0 \\
	a_{41} & a_{42} & 0 & 2a_{11} & 0 & 0 \\
	a_{51} & a_{52} & a_{42} & a_{32} & 3a_{11} & 0 \\
	a_{61} & a_{62} & -a_{32} & a_{42} & 0 & 3a_{11}
\end{pmatrix}.
\]
So, $I_{32}-I_{63} +I_{54} \in \gotha_0 \cap \gothu$, and for $I_{33}+I_{55} \in \Der(\frh_{10})$ we have
\[
[I_{32}-I_{63} +I_{54}, I_{33}+I_{55}] = -I_{32} - I_{63} - I_{54} \not\in \gotha_0 \cap \gothu.
\]

\subsubsection{The Lie algebra $\frh_{11}$}
The non-zero brackets of $\frh_{11}$ are
\[
[e_1, e_2] = e_4, \qquad [e_1, e_3]= e_5, \qquad [e_1,e_4] =[e_2,e_3] = e_6.
\]
We set $J_0(e_1) =e_1 +e_2, J_0(e_3)= e_4-e_3, J_0(e_5) =e_6$, and obtain $\gotha_0$ consisting of the matrices
\[
\begin{pmatrix}
	a_{11} & 0 & 0 & 0& 0& 0 \\
	0 & a_{11} & 0 & 0 & 0 & 0 \\
	a_{31} & a_{32} & 2a_{11} & 0 & 0 & 0 \\
	a_{41} & a_{42} & 0 &2a_{11} & 0 & 0 \\
	a_{51} & a_{52} & \tfrac{1}{2}a_{32}+a_{42}-a_{31} & a_{32} & 3a_{11} & 0 \\
	a_{61} & a_{62} &  -a_{31}-a_{32}-a_{41} & a_{42}-a_{31} & 0 & 3a_{11}
\end{pmatrix}.
\]
We have $I_{31}+I_{42}-I_{63} \in \gotha_0 \cap \gothu$ and $I_{21}+I_{65}\in \Der(\frh_{11})$. Computing their bracket yields
\[
[I_{31}+I_{42}-I_{63}, I_{21}+I_{65}] = I_{41} \not\in \gotha_0 \cap \gothu.
\]

\subsubsection{The Lie algebra $\frh_{12}$}
The non-zero brackets of $\frh_{12}$ are
\[
[e_1, e_2] = e_4, \qquad [e_1, e_3] = e_5, \qquad [e_2, e_4] = e_6.
\]
Consider the complex structure $J_0$ given by $J_0(e_1) = e_1 + e_2$, $J_0(e_3) = e_3 + 2e_4$,
$J_0(e_5) = e_6$. Then the stabiliser $\gotha_0 \subset \Der(\frh_{12})$ consists of matrices of
the form
\[
\begin{pmatrix}
	a_{11} & 0 & 0 & 0 & 0 & 0 \\
	0 & a_{11} & 0 & 0 & 0 & 0 \\
	a_{31} & a_{32} & 2a_{11} & 0 & 0 & 0 \\
	a_{41} & a_{42} & 0 & 2a_{11} & 0 & 0 \\
	a_{51} & a_{52} & -2a_{31} - 2a_{32} + a_{42} & a_{32} & 3a_{11} & 0 \\
	a_{61} & a_{62} & -2a_{32} & -a_{41} & 0 & 3a_{11}
\end{pmatrix}.
\]
The element $I_{31} + 2I_{42}$ is contained in $\gotha_0 \cap \gothu$ and $I_{33} + I_{55}$ is a
derivation of $\frh_{12}$. Computing their bracket in $\Der(\frh_{12})$ we get
\[
[I_{31} + 2I_{42}, I_{33} + I_{55}] = -I_{31} \not\in \gotha_0 \cap \gothu.
\]

\subsubsection{The Lie algebra $\frh_{13}$}
The non-zero brackets of $\frh_{13}$ are
\[
[e_1, e_2] = e_4, \qquad [e_1, e_3] = e_5, \qquad [e_1, e_4] = e_5, \qquad [e_2, e_4] = e_6.
\]
Let $J_0$ be the complex structure given by $J_0(e_1) = e_1 + 2e_2$, $J_0(e_3) = e_4$,
$J_0(e_5) = e_6$. The stabiliser $\gotha_0$ is comprised of matrices of the form
\[
\begin{pmatrix}
	a_{11} & 0 & 0 & 0 & 0 & 0 \\
	0 & a_{11} & 0 & 0& 0 & 0 \\
	a_{31} & a_{32} & 2a_{11} & 0 & 0 & 0 \\
	a_{41} & a_{42} &0 & 2a_{11} & 0 &0 \\
	a_{51} & a_{52} & a_{32} + a_{42} & a_{32}+ a_{42} & 3a_{11} & 0 \\
	a_{61} & a_{62} & -a_{31}-2a_{32} &-a_{41} & 0 & 3a_{11}
\end{pmatrix}.
\]
The element $2I_{31} - I_{32} + I_{42}$ is contained in $\gotha_0 \cap \gothu$ and
$I_{22} + I_{33} + I_{44} + I_{55} + 2I_{66}$ is a derivation of $\frh_{13}$. Computing their bracket
in $\Der(\frh_{13})$ we get
\[
[2I_{31} - I_{32} + I_{42}, I_{22} + I_{33} + I_{44} + I_{55} + 2I_{66}] = -2I_{31}
\not\in \gotha_0 \cap \gothu.
\]

\subsubsection{The Lie algebra $\frh_{14}$}
The non-zero brackets of $\frh_{14}$ are
\[
[e_1, e_2] = e_4, \qquad [e_1, e_3] = e_6, \qquad [e_1, e_4] = e_5, \qquad [e_2, e_4] = -e_6.
\]
Consider the complex structure $J_0$ given by $J_0(e_1) = e_2$, $J_0(e_3) = e_4$,
$J_0(e_5) = -2e_6$. The associated stabiliser $\gotha_0 \subset \Der(\frh_{14})$ consists of matrices of the form
\[
\begin{pmatrix}
	a_{11} & 0 & 0 & 0 & 0 & 0 \\
	0 & a_{11} & 0 &0&0& 0 \\
	a_{31} & a_{32} & 2a_{11} & 0 & 0 & 0 \\
	a_{41} & a_{42} & 0 & 2a_{11} & 0 & 0 \\
	a_{51} & a_{52} & 0& a_{42} & 3a_{11} & 0 \\
	a_{61} & a_{62} & a_{31} + a_{42} & a_{32} + a_{41} & 0 & 3a_{11}
\end{pmatrix}.
\]
The element $I_{32} - I_{41}$ is contained in $\gotha_0 \cap \gothu$ and $I_{21} - I_{65}$ is a
derivation of $\frh_{14}$. Computing their bracket in $\Der(\frh_{14})$ we get
\[
[I_{32} - I_{41}, I_{21} - I_{65}] = I_{31} \not\in \gotha_0 \cap \gothu.
\]

\subsubsection{The Lie algebra $\frh_{15}$}

The non-zero brackets of $\frh_{15}$ are
\[
[e_1,e_2]=e_4, \qquad [e_1,e_3]=e_5, \qquad [e_2,e_4]=-e_5, \qquad [e_1,e_4]=[e_2,e_3] = e_6.
\]
For $J_0(e_1) = e_2, J_0(e_3) = e_4, J_0(e_5) = e_6$, the stabiliser $\gotha_0$ consists of the matrices
\[
\begin{pmatrix}
	a_{11} & a_{12} & 0 &0 &0 & 0 \\
	-a_{12} & a_{11} & 0 & 0 & 0 & 0 \\
	a_{31} & a_{32} & 2a_{11} &0 & 0 & 0 \\
	a_{41} & a_{42} & 0 & 2a_{11} & 0 & 0 \\
	a_{51} & a_{52} & a_{42} - a_{31} & a_{32} + a_{41} & 3a_{11} & a_{12} \\
	-a_{52} & a_{51} & -a_{32} - a_{41} & a_{42} - a_{31} & -a_{12} & 3a_{11}
\end{pmatrix},
\]
We have $I_{31}- I_{53}- I_{64} \in \gotha_0 \cap \gothu$ and $I_{53} \in \Der(\frh_{15})$. Computing their brackets gives
\[
[I_{31}- I_{53}- I_{64}, I_{53}] = -I_{51} \not\in \gotha_0 \cap \gothu.
\]
So, $\gotha_0 \cap \gothu$ is not an ideal.

\subsubsection{The Lie algebra $\frh_{16}$}
The non-zero brackets of $\frh_{16}$ are
\[
[e_1, e_2] = e_4, \qquad [e_1, e_4]= e_5, \qquad [e_2,e_4] = e_6.
\]
For $J_0(e_1) =e_2, J_0(e_3) =e_4, J_0(e_5) =e_6$, the stabiliser $\gotha_0$ is comprised of matrices of the form
\[
\begin{pmatrix}
	a_{11} & a_{12} & 0 &0 & 0 & 0 \\
	-a_{12} & a_{11} & 0 &0 & 0 & 0 \\
	a_{31} & a_{32} & 2a_{11} & 0 &0 & 0 \\
	a_{41} & a_{42} & 0 & 2a_{11} & 0& 0 \\
	a_{51} & a_{52} &a_{32} & a_{42} & 3a_{11} &a_{12} \\
	-a_{52} & a_{51} &-a_{31} & -a_{41} & -a_{12} & 3a_{11}
\end{pmatrix} 
\]
Thus, we have $I_{31}-I_{63} \in \gotha_0 \cap \gothu$, and for $I_{31} \in \Der(\frh_{16})$ we compute
\[
[I_{31}-I_{63}, I_{31}] = -I_{61} \not\in \gotha_0 \cap \gothu.
\]

This case concludes the proof of the first item of Theorem \ref{thm: coarse moduli space dim 6}.
\subsection{The Lie algebras $\frh_{19}^-$ and $\frh_{26}^+$}

It follows from the classification of complex structures on $\frg \in \{ \frh_{19}^-, \frh_{26}^+\}$ proved in \cite{Ugarte14Non-nilpotent} that the $\Aut(\frg)$-orbits in $\kc(\frg)$ make up a union of connected components.

Therefore, we can write $\kt(\frg)$ as a union of connected components of the form $\Out(\frg)/H$ for some closed subgroup $H \subset \Out(\frg)$. Recall that $O_\Gamma$ denotes the image of the map $\Aut(\Gamma) \to \Out(\frg)$, and we have $\km(\frg, \Gamma) = \kt(\frg)/O_{\Gamma}$. Using Proposition~\ref{prop: transitive action}, we will show that $O_{\Gamma}$ acts properly discontinuously on $\kt(\frg)$.

\subsubsection{The Lie algebra $\frh_{19}^-$} For the Lie algebra $\frh_{19}^-$ it is convenient to work with a basis with the following structure equations:
\[
[e_1,e_3] = e_4, \quad [e_2,e_3] = e_5, \quad [e_1,e_4] = [e_2,e_5] = e_6.
\]
In the notation of \cite{Magnin07} this is the Lie algebra $M14_1$ and it was shown that the automorphism group $\Aut(\frh_{19}^-)$ is comprised of matrices of the form
\begin{equation*}
	\begin{pmatrix}
		a_{11} & \sigma a_{21} & 0&0 &0 &0 \\
		a_{21} & -\sigma a_{11} & 0&0 & 0 &0 \\
		0 & 0 & a_{33} & 0 & 0& 0\\
		a_{41} & a_{42} & 	a_{43}& a_{11}a_{33} & \sigma a_{21}a_{33} &0 \\
		\sigma(a_{21}u-a_{42}) & \sigma a_{41}-ua_{11}  & a_{53} & a_{21}a_{33}& -\sigma a_{11}a_{33} & 0 \\
		a_{61} & a_{62} & a_{63} & a_{64} & a_{65}  & (a_{11}^2 + a_{21}^2)a_{33}
	\end{pmatrix}
	\in \GL(6,\IR),
\end{equation*}
where $\sigma \in \{1,-1\}$, $u \in \IR$ and
\begin{align*}
	a_{64} &= a_{53}a_{21} + a_{43}a_{11},\\
	a_{65} &= -\sigma(a_{53}a_{11}- a_{43}a_{21}).
\end{align*}
Let $J_{\varepsilon}$ be the complex structure corresponding to the form
\[
\omega_\varepsilon = (e^1 +ie^2) \wedge (e^4+ie^5) \wedge (e^3+i\varepsilon e^6),
\]
where $\varepsilon \in \{1,-1\}$. Then the orbits $\Aut(\frh_{19}^-) \cdot J_\varepsilon$ form a union of connected components of $\kc(\frh_{19}^-)$, and the stabiliser of the $\Aut(\frh_{19}^-)$-action with respect to $J_\varepsilon$ is the group
\[
\widetilde{H} = \{\exp(\ad_{re_3}) \circ A(s,t) \, | \, r,t \in \IR, s \in \IR\backslash \{0\} \},
\]
where $\exp \colon \Der(\frg) \to \Aut(\frg)$ and
\[
A(s,t) =
\begin{pmatrix}
	\cos(t) & -\sin(t) & & & & \\
	\sin(t) & \cos(t) & & &  \textup{\huge 0}& \\
	&  & s &  & &\\
	&  & & s\cos(t)& -s\sin(t) &\\
	& \textup{\huge 0} &  & s\sin(t) & s\cos(t) &  \\
	&  &  &  & & s
\end{pmatrix}.
\]
Since $\exp(\ad_{e_3}$) is an inner automorphism, the image $H = \im(\widetilde H \to \Out(\frg))$ is the group $H = \{A(s,t)\, |\, s \in \IR\backslash \{0\}, t \in \IR\}$. Recall that $\Out^1(\frh_{19}^-)$ denotes the group of outer automorphisms whose determinant has absolute value $1$. So, we have $\Out^1(\frh_{19}^-) H = \Out(\frh_{19}^-)$ and the intersection $\Out^1(\frh_{19}^-) \cap H$ is compact. Therefore, the second item of Proposition~\ref{prop: transitive action}, implies that the group $O_\Gamma$ acts properly discontinuously on $\kt(\frh_{19}^-)$ for any lattice $\Gamma$.

\subsubsection{The Lie algebra $\frh_{26}^+$} The non-zero brackets of $\frh_{26}^+$ are
\[
[e_1,e_2] = e_3, \quad [e_1,e_3] = e_4, \quad [e_2,e_3] = e_5, \quad [e_1,e_4] =[e_2,e_5] = e_6.
\]
In the notation of \cite{Magnin07} this is the Lie algebra $M18_1$ and it was shown that the automorphism group $\Aut(\frh_{26}^+)$ is comprised of matrices of the form
\begin{equation*}
	\begin{pmatrix}
		a_{11} & \sigma a_{21} &0 & 0& 0&0 \\
		a_{21} & -\sigma a_{11} &0 &0 & 0 &0 \\
		a_{31} & a_{32} & -\sigma \delta &  0&0 &0\\
		a_{41} & a_{42} & 	a_{43}& -\sigma a_{11}\delta & -a_{21} \delta &0\\
		a_{51} & a_{52} & a_{53} & -\sigma a_{21} \delta & a_{11} \delta &  0\\
		a_{61} & a_{62} & a_{63} & a_{32}\delta & -a_{31} \delta & -\sigma \delta^2 
	\end{pmatrix}
	\in \GL(6,\IR),
\end{equation*}
where $\sigma \in \{1,-1\}$ and 
\begin{align*}
	\delta &= a_{11}^2 + a_{21}^2\\
	a_{43} &= a_{32}a_{11}- \sigma a_{31}a_{21},\\
	a_{53} &= a_{32}a_{21}+ \sigma a_{31}a_{11},\\
	a_{63} &= a_{42}a_{11}-\sigma a_{41}a_{21}+ \sigma a_{51}a_{11} + a_{52}a_{21}.
\end{align*}
As above in the previous case, we let $J_{\varepsilon}$ be the complex structure corresponding to the form
\[
\omega_\varepsilon = (e^1 +ie^2) \wedge (e^4+ie^5) \wedge (e^3+i\varepsilon e^6),
\]
where $\varepsilon \in \{1,-1\}$. Then the orbits $\Aut(\frh_{26}^+) \cdot J_\varepsilon$ form a union of connected components of $\kc(\frh_{26}^+)$, and the stabiliser of the $\Aut(\frh_{26}^+)$-action with respect to $J_{\varepsilon}$ is the group
	\[
	\widetilde{H} = \{\exp(\ad_{se_3}) \circ A(t) \, | \, s \in \IR\backslash\{0\},t \in \IR \},
	\]
	where
	\[
	A(t) =
	\begin{pmatrix}
		\cos(t) & -\sin(t) & & & & \\
		\sin(t) & \cos(t) & & &  \textup{\huge 0}& \\
		 &  & 1 &  & &\\
		 &  & & \cos(t)& -\sin(t) &\\
		 & \textup{\huge 0} &  & \sin(t) & \cos(t) &  \\
		 &  &  &  & & 1
	\end{pmatrix}.
	\]
Again, since $\exp(\ad_{e_3})$ is an inner automorphism, the image $H = \im(\widetilde H \to \Out(\frg))$ is compact, and thus by the first item of Proposition \ref{prop: transitive action}, $O_{\Gamma}$ acts properly discontinuously on $\kt(\frh_{26}^+)$ for every lattice.

This case concludes the proof of Theorem \ref{thm: coarse moduli space dim 6}.

\thispagestyle{fixedhead}
\printbibliography

\end{document}